\documentstyle {article}

\def\qmod#1#2{{\hbox{}^{\displaystyle{#1}}}\!\big/\!\hbox{}_{
\displaystyle{#2}}}

\newfam\msbfam

\font\tenmsb=msbm10
\font\sevenmsb=msbm10 at 7pt
\font\fivemsb=msbm10 at 5pt

\textfont\msbfam=\tenmsb
\scriptfont\msbfam=\sevenmsb
\scriptscriptfont\msbfam=\fivemsb

\def\Bbb{\fam\msbfam\tenmsb}

\def\C{{\Bbb C}}

\def\G{{\Bbb G}}
\def\H{{\Bbb H}}

\def\P{{\Bbb P}}

\def\R{{\Bbb R}}
\def\Z{{\Bbb Z}}

\def\union{\mathop{\bigcup}}
\def\qed {\hfill\vrule height6pt width6pt depth0pt \bigskip}

\def\map{\longrightarrow}
\def\textmap#1{\mathop{\vbox{\ialign{
                                ##\crcr
    ${\scriptstyle\hfil\;\;#1\;\;\hfil}$\crcr
    \noalign{\kern-1pt\nointerlineskip}
    \rightarrowfill\crcr}}\;}}

\def\textlmap#1{\mathop{\vbox{\ialign{
                                ##\crcr
    ${\scriptstyle\hfil\;\;#1\;\;\hfil}$\crcr
    \noalign{\kern-1pt\nointerlineskip}
    \leftarrowfill\crcr}}\;}}

\newfam\meuffam

\font\tenmeuf=eufm10
\font\sevenmeuf=eufm7

\textfont\meuffam=\tenmeuf
\scriptfont\meuffam=\tenmeuf
\scriptscriptfont\meuffam=\sevenmeuf

\def\germ{\fam\meuffam\tenmeuf}

\def\cg{{\germ c}}

\def\g{{\germ g}}
\def\hg{{\germ h}}

\def\tg{{\germ t}}

\begin{document}
\def\Pr{{\rm Pr}}
\def\tr{{\rm Tr}}
\def\End{{\rm End}}
\def\Spin{{\rm Spin}}
\def\U{{\rm U}}
\def\SU{{\rm SU}}
\def\SO{{\rm SO}}
\def\PU{{\rm PU}}
\def\GL{{\rm GL}}
\def\spin{{\rm spin}}
\def\u{{\rm u}}
\def\su{{\rm su}}
\def\so{{\rm so}}
\def\pu{{\rm pu}}
\def\Pic{{\rm Pic}}
\def\Iso{{\rm Iso}}
\def\NS{{\rm NS}}
\def\deg{{\rm deg}}
\def\Hom{{\rm Hom}}
\def\h{{\germ h}}
\def\Herm{{\rm Herm}}
\def\Vol{{\rm Vol}}
\def\pf{{\bf Proof: }}
\def\id{{\rm id}}
\def\i{{\germ i}}
\def\im{{\rm im}}
\def\rk{{\rm rk}}
\def\ad{{\rm ad}}
\def\h{{\bf H}}
\def\coker{{\rm coker}}
\def\dv{\bar\partial}
\def\Ad{{\rm Ad}}
\def\ad{{\rm ad}}
\def\dva{\bar\partial_A}
\def\da{\partial_A}
\def\p{\partial\bar\partial}
\def\pa{\partial_A\bar\partial_A}
 \def\Dr{{\raisebox{.17ex}{$\not$}}\hskip -0.4mm{D}}
\def\gr{{\scriptscriptstyle|}\hskip -4pt{\g}}
\def\subsetint{{\  {\subset}\hskip -2.45mm{\raisebox{.28ex}
{$\scriptscriptstyle\subset$}}\ }}
\def\nr{\parallel}
\newtheorem{sz}{Satz}[subsection]
\newtheorem{thry}[sz]{Theorem}
\newtheorem{pr}[sz]{Proposition}
\newtheorem{re}[sz]{Remark}
\newtheorem{co}[sz]{Corollary}
\newtheorem{dt}[sz]{Definition}
\newtheorem{lm}[sz]{Lemma}
\newtheorem{cl}[sz]{Claim}

\title{Non-Abelian Seiberg-Witten Theory and stable oriented pairs}
\author{ Andrei Teleman\thanks{Partially supported  by: AGE-Algebraic
Geometry in Europe,
contract No ERBCHRXCT940557 (BBW 93.0187), and by  SNF, nr. 21-36111.92}  }
\date{ }
\maketitle

\setcounter{section}{-1}
\section{Introduction}

The aim of this paper is to develop a systematic  theory  of non-abelian
Seiberg-Witten
equations. The   equations we introduce and study are associated with a
$Spin^G(4)$-structure on a 4-manifold, where $G$ is a closed subgroup of the
unitary group
$U(V)$ containing the central involution $-\id_V$. We call these equations the
$G$-monopole
equations. For $G=S^1$, one recovers the classical (abelian) Seiberg-Witten
equations
[W], and
the  case $G=Sp(1)$ corresponds to the "quaternionic monopole equations"
 introduced in [OT5].
Fixing the determinant of the connection component in the $U(2)$-monopole
 equations, one
gets the so called $PU(2)$-monopole equations, which should be regarded
as a twisted version
of quaternionic monopole equations and will be extensively studied  in the
second part of this
paper.

It is known ([OT4], [OT5], [PT2]) that the most natural way to prove the
equivalence between
Donaldson theory and Seiberg-Witten theory is to consider a suitable
moduli space of
non-abelian monopoles. In [OT5]  it was shown that an $S^1$-quotient of
 a moduli space of
quaternionic monopoles should give an homological equivalence between
a fibration over  a
union of Seiberg-Witten moduli spaces  and  a fibration over   certain
$Spin^c$-moduli spaces
  [PT1].

By the same method, but using moduli spaces of $PU(2)$-monopoles
 instead of quaternionic
monopoles,  one should be able to express any Donaldson polynomial
 invariant
  in terms of Seiberg-Witten invariants associated with the
\underbar{twisted} abelian
monopole equations of [OT6].

The idea can be extended to get information about the
Donaldson theories associated with an arbitrary symmetry
group $G$, by relating the
corresponding polynomial invariants to Seiberg-Witten-type
invariants associated with smaller
symmetry groups. One has only to consider a suitable moduli
 space of $G$-monopoles and to
notice that this moduli space contains distinguished closed
 subspaces of "reducible solutions".
The reducible solutions with trivial spinor-component can
be identified with $G$-instantons,
and all the others reductions can be regarded as  monopoles
associated to a smaller group.

It is important to point out that, if the base manifold is a
K\"ahler surface one has
Kobayashi-Hitchin-type correspondences (see  [D], [DK], [K],
 [LT] for the instanton case) which
give  a complex geometric description of the moduli spaces
 of $SU(2)$, $U(2)$ or
$PU(2)$-monopoles (see section 2). The first two cases were
already studied in  [OT5] and
[OT1]. In the algebraic case    one can explicitly compute such
moduli spaces  of non-abelian
monopoles and prove   the existence of a projective
compactification. The points
corresponding to instantons and abelian monopoles can be
easily identified (see also [OST]).\\

The theory has interesting extensions to manifolds of other dimensions.
 On Riemann surfaces
for instance, one can use moduli spaces of $PU(2)$-monopoles to reduce
   the
computation of the volume or   the Chern numbers of a moduli space of
 semistable rank 2-
bundles to computations on the symmetric powers of the base, which
 occur in the moduli
space of $PU(2)$-monopoles as  subspaces of abelian reductions. \\

The present paper is divided into two parts: The first   deals with the
general theory of
$Spin^G$-structures and $G$-monopole equations. We give classification
 theorems for
$Spin^G$-structures in principal bundles, and an   explicit description of
the set of equivalence
classes in the cases $G=SU(2)$, $U(2)$, $PU(2)$. Afterwards we introduce
the $G$-monopole
equations in a natural way  by coupling the Dirac harmonicity condition for
 a pair formed by a
connection and a spinor, with the vanishing condition for a
generalized moment map. This first
part ends with a section dedicated to the concept of reducible
solutions of   the
$G$-monopole equations. Describing the  moduli spaces of
 $G$-monopoles around the
reducible loci  is the first step   in order to  express the
 Donaldson
invariants associated with the   symmetry group $G$ in terms
of Seiberg-Witten-type
reductions.\\

In the second part of the paper, we give a complex geometric
interpretation of the moduli
spaces of $PU(2)$-monopoles in terms of   stable oriented pairs,
by proving a
Kobayashi-Hitchin type correspondence. Using this result,  we describe
 a simple
example of moduli space of $PU(2)$-monopoles on $\P^2$, which
illustrates in a concrete
case how our moduli spaces  can be used to relate    Donaldson and
  Seiberg-Witten
invariants.

In order to be able to give general explicit  formulas relating    the
Donaldson polynomial
invariants to  Seiberg-Witten invariants, it remains   to construct
$S^1$-equivariant smooth  perturbations of the moduli spaces of
$PU(2)$-monopoles, to
construct an Uhlenbeck compactification of the perturbed moduli spaces,
 and finally to
give explicit descriptions of the ends of the (perturbed) moduli spaces.

The first two problems are treated in [T1], [T2]. Note that the proof of the
corresponding
transversality results for other moduli spaces of non-abelian connections
coupled with
harmonic spinors ([PT1], [PT2]) are not complete ([T1]).  The third problem, as
well as
generalizations to larger symmetry groups will be treated in a future paper.

 I  thank Prof. Christian Okonek for encouraging me to write this paper,
as well
as for the careful reading of the text and his valuable suggestions.

\section{G-Monopoles on 4-manifolds}

\subsection{The group $Spin^G$ and $Spin^G$-structures}

\subsubsection{$Spin^G$-structures in principal bundles}

Let $G\subset U(V)$ be a closed subgroup of the unitary group of
a Hermitian vector space $V$,
 suppose that $G$ contains the central involution $-\id_V$, and denote
 by $\g\subset u(V)$ the
Lie algebra of $G$. We put
$$Spin^G :=Spin \times_{\Z_2} G \ .
$$

By   definition we get the following fundamental   exact sequences:
$$\begin{array}{c}1\map  Spin \map  Spin^G \stackrel{\delta}\map
\qmod{G}{\Z_2}\map
1
\\
 1\map G\map  Spin^G \stackrel{\pi}{\map } SO \map  1
\\
 1\map  \Z_2\map  Spin^G\textmap{(\pi,\delta)} SO\times
\qmod{G}{\Z_2}\map  1
\end{array}\eqno{(*)}$$

Note first that there are well defined  morphisms
$$\ad_G:Spin^G\map O(\g)\  ,\  \Ad_{G}:Spin^G \map {\rm Aut}(G)
$$
induced by the morphisms $\ad:G\map O(\g)$ and $\Ad:G\map
{\rm Aut}( G)$.

If $P^G$ is principal $Spin^G$-bundle, we denote by $\G(P^G)$,
${\scriptscriptstyle|}\hskip
-4pt{\g}(P^G)$     the fibre bundles
$P^G\times_{\Ad_G}G$,
$P^G\times_{\ad_G}\g$.  The group of sections
$${\cal G}(P^G):=\Gamma(\G(P^G))
$$
in $\G(P^G)$ can be identified with the group of bundle-automorphisms
of $P^G$ over
the $SO$-bundle $P^G\times_\pi SO$. After a suitable Sobolev completion
${\cal
G}(P^G)$
becomes a Hilbert Lie group, whose Lie algebra is the corresponding Sobolev
completion
of
$\Gamma(\gr(P^G))$.

We put also
$$ \delta(P^G):=P^G\times_{\delta} \left(\qmod{G}{\Z_2}\right)\ .$$
Note that $\G(P^G)$ can be identified with the bundle
$\delta(P^G)\times_{\bar\Ad} G$
associated with the $\qmod{G}{\Z_2}$-bundle $\delta(P^G)$.

Let $P$ a principal $SO(n)$-bundle over a topological space $X$. A
\underbar{$Spin^G(n)$}-\underbar{structure}  in $P$ is a bundle
morphisms $P^G\map  P$ of
type $\pi$, where $P^G$ is a principal $Spin^G(n)$-bundle over $X$.
Equivalently, a
$Spin^G(n)$-structure in
$P$ can be regarded as a pair consisting  of a
$Spin^G(n)$-bundle $P^G$ and an orientation preserving linear
isometry
$$\gamma:P\times_{SO(n)}\R^n\map P^G\times_{\pi}\R^n $$
(called the \underbar{Clifford} \underbar{map} of the
structure).

Two $Spin^G$-structures $P^G_0\textmap{\sigma_0} P$,
$P_1^G\textmap{\sigma_1} P$ in $P$ are called \underbar{equivalent},
if the $Spin^G$
bundles
$P^G_0$, $P_1^G$ are isomorphic   over    $P$.

If $(X,g)$ is an oriented  Riemannian  $n$-manifold, a
$Spin^G(n)$-structure in $X$ is a $Spin^G(n)$-structure $P^G\map  P_g$ in
the bundle $P_g$ of oriented $g$-orthonormal coframes of $ X $.   This is
equivalent with
the data of a pair $(P^G,\gamma)$, where $P^G$ is a $Spin^G(n)$-bundle and
$\gamma:\Lambda^1_X\stackrel{\simeq}{\map} P^G\times_\pi\R^n$ is a linear
orientation-preserving isometry. Here
$\Lambda^1_X$ stands for the cotangent bundle of $X$, endowed with the dual
$SO(n)$-structure.
\vspace{5mm}

Let $X$ be a fixed paracompact topological space. Note that there is a
natural map
$H^1(X,\underline {{G}/{\Z_2}}) {\map} H^2(X,\Z_2)$, which we denote by $w$. If
$G=Spin(k)$,
$w$ coincides with the usual morphism $w_2$ defined on the set of
$SO(k)$-bundles.
By the
third exact sequence in
$(*)$ we get the following simple classification result
\begin{pr}  The map $P^G\longmapsto
(P^G\times_\pi SO,\delta(P^G))$ defines a surjection of the set of
isomorphism classes of
$Spin^G$-bundles onto the set of isomorphism classes of pairs $(P,\Delta)$
consisting   of an
$SO$-bundle and a $\qmod{G}{\Z_2}$-bundle satisfying $w_2(P)+w(\Delta)=0$.
Two
$Spin^G$-bundles have the same image if and only if they are congruent modulo
the natural
action of $H^1(X,\Z_2)$ in $H^1(X,\underline{Spin^G})$.
\end{pr}
\pf Indeed, the natural morphism
$H^1(X,\underline{SO_{}}\times\underline{G/\Z_2})
\map
H^2(X,\Z_2)$ is given by $(P,\Delta)\longmapsto (w_2(P)+w(\Delta))$.
\qed
\\
 For instance, we have the following result
\begin{pr} Let $X$ be a 4-manifold. The group $H^1(X,\Z_2)$ acts trivially
on the
 set of
(equivalence classes of) $Spin^c(4)$-bundles over $X$. Equivalence classes of
$Spin^c(4)$-bundles over $X$   are classified by pairs $(P,\Delta)$
consisting of
 an
$SO(4)$-bundle $P$ and an
$S^1$-bundle
$\Delta$ with $w_2(P)+w_2(\Delta)=0$.
\end{pr}
\pf Using the identification (see [OT1], [OT3])
$$Spin^c(4)=\{(a,b)\in U(2)\times U(2)|\ \det a=\det b\}\ ,$$
we get an
exact sequence
$$1\map Spin^c(4)\map U(2)\times U(2)\map  S^1\map  1\ .
$$
Using this, one can prove that, on 4-manifolds, the  data of an
(equivalence class of)
$Spin^c(4)$-bundles
 is equivalent
to the data of a pair of $U(2)$-bundles having isomorphic determinant line
bundles. The
action of $H^1(X,\Z_2)$ is given by tensoring with flat line bundles with
structure group
$\Z_2$. The Chern class of such line bundles is 2-torsion, hence the assertion
 follows
from the classification of unitary vector bundles on 4-manifolds in terms
of Chern
classes.
\qed

The classification of the $Spin^G$-structures in a given $SO$-bundle $P$ is a
more delicate
problem.
\begin{pr} Fix a $Spin^G$-structure $\sigma:P^G\map P$ in $P$. Then
the set of equivalence classes of $Spin^G$-structures in $P$ can be identified
with the
cohomology set $H^1(X,\G(P^G ))$ of the sheaf of sections in the bundle
$\G(P^G )$.
\end{pr}

Recall that $\G(P^G)$ can be identified with the bundle $\delta(P^G)
\times_{\bar\Ad} G$
associated with the $\qmod{G}{\Z_2}$-bundle $\delta(P^G)$.

Therefore we get the exact sequence of   bundles of groups
$$1\map \Z_2\map\G(P^G)\map \delta(P^G)\times_{\Ad}
\left(\qmod{G}{\Z_2}\right)\map 1\ .
$$
The third term coincides with the gauge group of automorphisms of
$\delta(P^G)$. The cohomology set
$H^1\left(X,\delta(P^G)\times_{\bar\Ad}
\left(\qmod{G}{\Z_2}\right)\right)$ of the
associated sheaf coincides with the pointed set of (equivalence classes of)
$\qmod{G}{\Z_2}$-bundles over
$X$ with distinguished element $\delta(P^G)$. This shows that
$\qmod{H^1(X,\G(P^G))}{H^1(X,\Z_2)}$ can be identified with the set of
$\qmod{G}{\Z_2}$-bundles $\Delta$ with $w(\Delta)=w(\delta(P^G))$.
Therefore
\begin{pr}  The map
$$(\sigma:P^G\map P) \longmapsto \delta(P^G)$$
is a surjection of the set of (equivalence classes of) $Spin^G$-structures
in $P$ onto the
set of
$\qmod{G}{\Z_2}$-bundles $\Delta$ satisfying  $w(\Delta)+w_2(P)=0$.
Two
$Spin^G$-structures have the same image if and only if they are congruent
 modulo the
natural action of
$H^1(X,\Z_2)$.
\end{pr}
Proposition 1.1.2 and the proposition below show  that the
 classification of
$Spin^G$-structures in the
$SO$-bundle $P$ is  in general different from the classification
 of $Spin^G$-bundles with
associated $SO$-bundle  isomorphic to $P$.
\begin{pr}\hfill{\break}
1. If $G=S^1$ then $Spin^{S^1}=Spin^c$,   $\G(P^G)=X\times S^1$,
hence the
set of $Spin^c$-structures in $P$ is a $H^1(X,\underline{S}^1)=
H^2(X,\Z)$-torsor if it is
non-empty. The $H^1(X,\Z_2)$-action in the set of $Spin^c$-structures
 in $P$, factorizes
through a ${\rm Tors}_2 H^2(X,\Z)$-action, which is free and whose
 orbits  coincide with
the fibres of the determinant map $(\sigma:P^c\map P)\longmapsto
\delta(P^c)$.\\
 2. Suppose that $X$ is a 4-manifold, $P$ is an $SO$-bundle over $X$
and that $G$ is one of
the following:\\ a) $SU(r)$, $r\geq 2$,
b) $U(r)$, $r\geq 2$, $r$ even.
c) $Sp(r)$, $r\geq 1$.

Then $H^1(X,\Z_2)$ acts trivially in the set of $Spin^G$-structures in
 $P$, hence the
classification of  $Spin^G$-structures in $P$ reduces to the classification of
$\qmod{G}{\Z_2}$-bundles over $X$.
\end{pr}

\pf \\
1. The first assertion follows immediately from Propositions 1.1.3 and 1.1.4. \\
  2. Let
$\sigma_i:P^G_i\map  P$, $i=0,\ 1$ be two $Spin^G$-structures in $P$.
We consider the
locally trivial bundle $Iso_P(P^G_1,P^G_0)$ whose fibre in $x\in X$ consists
of isomorphism
$\rho_x:(P^G_1)_x\map (P^G_0)_x$ of right $Spin^G$-spaces  which make the
 following
diagram commutative.
$$\begin{array}{rcl}
(P^G_1)_x&\stackrel{\rho_x}{\map }&{(P^G_0)_x}_{\phantom{X_{X_{X_X}}}  }\\
{\scriptstyle\sigma_{0x}}\searrow&&\swarrow{\scriptstyle
\sigma_{\scriptscriptstyle
1x}}\\
& P_x&
\end{array}
$$
$Iso_P(P^G_1,P^G_0)$ is a principal bundle in the sense of Grothendieck
with structure
group bundle $\G(P^G_0)$. The $Spin^G$-structures $\sigma_i$ are equivalent
if and only if
$Iso_P(P^G_1,P^G_0)$ admits a section.  Consider first the case $G=SU(r)$
($r\geq 2$) or
$Sp(r)$ ($r\geq 1$). Since
$\pi_i\left([Iso_P(P^G_1,P^G_0)]_x\right)=0$ for $i\leq 2$ and
$\pi_3\left([Iso_P(P^G_1,P^G_0)]_x\right)$ can be canonically identified
with $\Z$,
the
obstruction $o(\sigma_1,\sigma_0)$ to the existence of such a section is an
element in
$H^4(X,\Z)$. Assume now that $\sigma_1=\lambda\sigma_0$ for some
$\lambda\in
H^1(X,\Z_2)$ and let $p:\tilde X\map  X$ the cover associated to
$\ker\lambda\subset\pi_1(X)$. It is easy to see that one has
$o(p^*(\sigma_1),p^*(\sigma_0))=p^*(o(\sigma_1,\sigma_0))$. But, since
$p^*(\lambda)=0$,  we get $p^*(\sigma_1)=p^*(\sigma_0)$ hence
$o(p^*(\sigma_1),p^*(\sigma_0))=0$. Since $p^*:H^4(X,\Z)\map
H^4(\tilde X,\Z)$ is injective for a 4-manifold $X$, the assertion follows
 immediately.

Finally consider $G=U(r)$.  When $r\geq 2$ is even,  the determinant map
 $U(r)\map S^1$
induces a morphism $Spin^{U(r)}
\map S^1$. If  $\sigma_1=\lambda\sigma_0$, then there is a natural
 identification
$P^G_1\times_{\det} S^1=P^G_0\times_{\det} S^1$, hence, denoting this
line bundle by $L$,
we get a subbundle
$Iso_{P,L}(P^G_1,P^G_0)$ of $Iso_P(P^G_1,P^G_0)$ consisting fibrewise of
isomorphisms
$(P^G_1)_x\map (P^G_0)_x$ over $P_x\times L_x$. Since the standard fibre of
$Iso_{P,L}(P^G_1,P^G_0)$ is $SU(r)$, the same argument as above shows that
 this bundle
admits a section, hence $\sigma_1$ and $\sigma_0$ are equivalent.

\qed

\subsubsection{$Spin^G(4)$-structures on 4-manifolds and spinor bundles}

Let $\H_{\pm}$ be two copies of the quaternionic skewfield, regarded as right
 quaternionic
vector spaces.    The canonical left actions of $Sp(1)$ in $\H_{\pm}$ define an
orthogonal representation of the group
$$Spin(4)=Sp(1)\times Sp(1)=SU(2)\times SU(2)$$
in $\H\simeq \Hom_{\H}(\H_+,\H_-)$, which gives the standard identification
$$\qmod{SU(2)\times SU(2)}{\Z_2}=SO(\H)=SO(4)$$

Therefore, the group $Spin^G(4)=\qmod{SU(2)\times SU(2)\times G}{\Z_2}$
comes with  2
unitary representations
$$\lambda_\pm: Spin^G(4)\map U(\H_{\pm}\otimes_\C V)$$
obtained by coupling the natural representation  of $G$ in $V$ with the
spinorial
representations $p_{\pm}:Spin(4)=SU(2)\times SU(2)\map  SU(2)$.

There are well defined adjoint morphisms
$$ {\rm ad}_{\pm}:Spin^G(4)\map O(su(2))\ ,\ \ \Ad_{\pm}:Spin^G(4)\map
{\rm Aut}(SU(2))$$
induced by the projections $p_{\pm}$ and the corresponding adjoint
 representations
associated with the Lie group $SU(2)$. If $P^G$ is a $Spin^G(4)$-bundle,
we denote by
$\ad_{\pm}(P^G)$,   $\Ad_{\pm}(P^G)$ the corresponding  bundles with
fibres $su(2)$,
$SU(2)$ associated with $P^G$.
The \underbar{spinor} \underbar{vector} \underbar{bundles} associated
 with a
$Spin^G(4)$-bundle $P^G$ are defined by
$$\Sigma^{\pm}=\Sigma^{\pm}(P^G):=P^G\times_{\lambda_{\pm}}
 (\H_{\pm}\otimes_\C
V)\ ,
 $$

The bundles $\ad_{\pm}(P^G)$, $\gr(P^G)$ are real subbundles of the
endomorphism bundle
$\End_\C(\Sigma^{\pm})$. The bundle $\G(P^G)$ acts fibrewise unitarily
 in the
bundles $\Sigma^{\pm}$. On the other hand, the identification $\H\simeq
\Hom_\H(\H_+,\H_-)$ defines a real embedding
 \begin{equation}
P^G\times_\pi\H\stackrel{ }{\longrightarrow}\Hom_{\G(P^G)}(\Sigma^+,
\Sigma^-)\subset\Hom_\C(
\Sigma^+,\Sigma^-)
\end{equation}
of the $SO(4)$-vector bundle $P^G\times_\pi\H$ in the bundle
$\Hom_{\G(P^G)}(\Sigma^+,\Sigma^-)$ of
$\C$-linear morphisms $\Sigma^+\map \Sigma^-$ which commute with the
$\G(P^G)$-action.

The data of a $Spin^G(4)$-structure with principal bundle $P^G$ on an
oriented Riemannian 4-manifold $X$ is   equivalent to the data  of an
 orientation-preserving
isomorphism $\Lambda^1_X\stackrel{\gamma}{\map}P^G\times_\pi\H$,
which defines (via
the monomorphism in (1)) a
\underbar{Clifford}
\underbar{multiplication}
$(\Lambda^1\otimes\C)\otimes \Sigma^+
\map \Sigma^-$ commuting with the ${\G(P^G)}$ actions in $\Sigma^{\pm}$.
Moreover,
as in the classical $Spin^c(4)$ case [OT1], [OT3], we also have induced
identifications (which multiply the norms by 2)
$$\Gamma:\Lambda^2_{\pm}\map  \ad_{\pm}(P^G)\ .
$$
\subsubsection{ Examples}
1.  $Spin^c(4)$-structures:\\

The group $Spin^c(4):=Spin^{U(1)}(4)$  can be identified with the
subgroup
$$G_2:=\{(a,b)\in U(2)\times U(2)|\ \det a=\det b\}
$$
of $U(2)\times U(2)$. Via this identification, the map $\delta:Spin^c(4)\map
S^1\simeq\qmod{S^1}{\Z_2}$ in the exact sequence
$$1\map  Spin(4)\map  Spin^c(4)\stackrel{\delta}{\map} S^1\map  1
$$
is given by the formula $\delta(a,b)=\det a=\det b$.  The spinor bundles come
with
identifications
$$\det\Sigma^+\stackrel{\simeq}{\rightarrow}\det\Sigma^-
\stackrel{\simeq}{\rightarrow}
P^c\times_\delta\C\ .$$

The $SO(4)$-vector bundle $P^c\times_\pi \H$ associated with a
$Spin^c(4)$-bundle
$P^c$ can be identified with the bundle
$\R SU(\Sigma^+\Sigma^-)\subset\Hom(\Sigma^+,\Sigma^-)$  of   real
multiples of isometries of determinant 1.

Using these facts, it easy to see that a $Spin^c(4)$-structure can be recovered
from the
data of the spinor bundles, the identification between the determinant line
bundles
and the
Clifford map. More precisely
\begin{pr} The data  of a  $Spin^c(4)$-structure  in the
$SO(4)$-bundle $P$ over $X$ is equivalent to the data of a
triple consisting of:\\
i) A pair of $U(2)$-vector bundles $\Sigma^{\pm}$.\\
ii)  A unitary  isomorphism $\det\Sigma^+\stackrel{\iota}{\rightarrow}\det
\Sigma^-$.\\
iii) An  orientation-preserving linear isometry
$$\gamma:P\times_{SO(4)}\R^4 \rightarrow\R SU(\Sigma^+,\Sigma^-)\ .$$
\end{pr}
%
\pf  Given a triple $(\Sigma^{\pm},\iota,\gamma)$, we define $P^c$ to be the
manifold
over $X$
$$
\begin{array}{cl}
P^c:=\left\{ [x,(e_1^+,e_2^+),(e_1^-,e_2^-) ]|\right.&  x\in X,\
(e_1^\pm,e_2^\pm)\ {\rm  an\ orthonormal\ basis\ in\ }\Sigma^{\pm}_x,\\
&\left. \iota_*( e_1^+\wedge e_2^+)=e_1^-\wedge e_2^-\right\}\ .
\end{array}
$$

Every   triple $[x,(e_1^+,e_2^+),(e_1^-,e_2^-)]\in P^c_x$ defines an
orthonormal orientation-compatible basis in  $\R SU(\Sigma^+_x,\Sigma^-_x)$
which is
given  with respect to the frames $(e_1^\pm,e_2^\pm)$  by the Pauli matrices.
 Using
 the isomorphism $\gamma$, we get a bundle morphism from $P^c$ onto the
orthonormal
oriented frame bundle of
$P\times_{SO(4)}\R^4$, which can be canonically identified with $P$.
\qed

Let $P$ be a principal  $SO(4)$-bundle,   $P^c\stackrel{\cg_0} \map P$
a fixed
$Spin^c(4)$-structure in   $P$, $\Sigma^{\pm}$ the associated spinor
bundles, and
$$\gamma_0:P\times_{SO(4)}\R^4\map  P^c\times_\pi\H=
\R SU(\Sigma^+,\Sigma^-)$$
 the
corresponding Clifford map. For every $m\in H^2(X,\Z)$ let $L_m$ be
a Hermitian line
bundle of Chern class
$m$.  The fixed identification $\det \Sigma^+\textmap{\simeq}\det \Sigma^-$
 induces
an identification   $\det \Sigma^+\otimes L_m\textmap{\simeq}\det
\Sigma^-\otimes
L_m$, and the map
$$\gamma_m: P\times_{SO(4)}\R^4\map \R SU(\Sigma^+\otimes L_m,\Sigma^-
\otimes L_m)\ ,\
 \gamma_m(\eta):=\gamma_0(\eta)\otimes\id_{L_m}$$
 is the Clifford map of a
$Spin^c(4)$-structure $\cg_m$ in $P$ whose spinor bundles are $\Sigma^{\pm}
\otimes L_m$.

Using the results in the previous section (see also  [H], [OT1], [OT6])  we
get

%
\begin{pr} \hfill{\break}
i) An $SO(4)$-bundle $P$ admits a $Spin^c(4)$-structure iff $w_2(P)$ admits
 integral
lifts.\\
ii)   The set of isomorphism classes of $Spin^c(4)$-structures in an $SO(4)$-
bundle $P$ is
either empty or is an $H^2(X,\Z)$-torsor. If
$\gamma_0$ is a fixed
$Spin^c(4)$-structure in the
$SO(4)$-bundle
$P$, then the map
$m\longmapsto \cg_m$ defines a bijection between $H^2(X,\Z)$ and
the set  of
(equivalence classes of) $Spin^c(4)$-structures in $P$. \\
iii)   [HH] If $(X,g)$ is a compact oriented Riemannian 4-manifold,
then
$w_2(P_g)$ admits  integral lifts. In particular any compact oriented
Riemannian
4-manifold admits $Spin^c(4)$-structures\\
\end{pr}
2. $Spin^h(4)$-structures: \\

The quaternionic spin group is defined  by $Spin^h:=Spin^{Sp(1)}$.
By the classification results 1.1.4., 1.1.5 we get

\begin{pr} Let $P$ be an $SO$-bundle over a compact oriented 4-manifold $X$.
The map
$$\left[\sigma:P^h\map  P\right]\longmapsto [\delta(P^h)]$$
defines a 1-1 correspondence between the set of  isomorphism classes of
$Spin^h$-structures in $P$ and the set of isomorphism classes of
$PU(2)$-bundles $\bar
P$ over $X$  with $w_2(\bar P)=w_2(P)$.  The latter set can be
identified ([DK], p.41) with
$$\{p\in\Z|\ p\equiv w_2(P)^2\ {\rm mod}\ 4\}$$
via the Pontrjagin class-map.
 \end{pr}

In dimension 4, the group $Spin^h(4)$ can be identified with the
quotient
$$\qmod{SU(2)\times SU(2) \times SU(2)}{\{\pm(\id,\id,\id)\}}\ ,$$
 hence there is an exact sequence
\begin{equation}1\map \Z_2\map  SU(2)\times SU(2) \times SU(2)\map
Spin^h(4)\map
1\ .
\end{equation}

Let $G_3$ be the group
$$G_3:=\{(a,b,c)\in U(2)\times U(2)\times U(2)|\ \det a=\det b=\det c\}
$$
We have an exact sequence
$$1\map  S^1\map  G_3\map Spin^h(4)\map 1$$
extending the exact sequence (2). If $X$ is any manifold, the induced
map
$H^1(X,\underline{Spin^h(4)})\map H^2(X,\underline{\phantom{(}S^1})=
H^3(X,\Z)$
factorizes as
$$H^1(X,\underline{Spin^h(4)})\stackrel{\pi}{\map}
 H^1(X,\underline{SO(4)})\stackrel{w_2}{\map} H^2(X,\Z_2)\map
 H^2(X,S^1)\ .$$

Therefore a $Spin^h(4)$-bundle $P^h$ admits an $G_3$-reduction iff
the second
Stiefel-Whitney class $w_2(P^h\times_\pi SO(4))$ of the associated
$SO(4)$-bundle
admits an integral lift. On the other hand, the data of a $G_3$-structure
in a
$SO(4)$-bundle $P$ is equivalent to the data of a triple consisting of  a
$Spin^c(4)$-structure
$P^c\map P$ in $P$, a $U(2)$-bundle $E$, and an isomorphism
$$P^c\times_\delta\C\textmap{\simeq}\det E\ .$$
Therefore ( see [OT5]),

\begin{pr} Let $P$ be  a principal $SO(4)$-bundle whose second Stiefel-Whitney
 class $w_2(P)$
admits an integral lift. There is a 1-1 correspondence between isomorphism
classes of
$Spin^h(4)$-structures in $P$ and   equivalence classes of triples consisting
of  a
$Spin^c(4)$-structure $P^c\map P$ in $P$, a $U(2)$-bundle $E$, and an
isomorphism
$P^c\times_\delta \C\textmap{\simeq}\det E$.

Two   triples are equivalent  if,
after tensoring the first  with an $S^1$-bundle, they become isomorphic
over $P$.
\end{pr}

Let us identify\  $\qmod{SU(2) \times SU(2)}{\Z_2}$ with $SO(4)=SO(\H)$ as
explained
 above, and
denote by
$$\pi_{ij}:Spin^h\map SO(4) \ \ \ 1\leq i<j\leq 3$$
the three epimorphisms associated with the three projections of the product
$SU(2)\times SU(2)\times SU(2)$ onto $SU(2)\times SU(2)$. Note that
$\pi_{12}=\pi$.
The spinor bundles
$\Sigma^{\pm}(P^h)$ associated with a principal $Spin^h(4)$-bundle $P^h$ are
$$\Sigma^+(P^h)=P^h\times_{\pi_{13}}\C^4\ ,\  \   \Sigma^-(P^h)=
P^h\times_{\pi_{23}}\C^4 $$

This shows in particular that the Hermitian 4-bundles $\Sigma^{\pm}(P^h)$
come with
\underbar{a} \underbar{real} \underbar{structure} and   compatible
trivializations of
$\det(\Sigma^{\pm}(P^h))$.

Suppose now that the \ $Spin^h(4)$-bundle $P^h$ \ admits a  $G_3$-lifting\
, consider
the associated triple
$(P^c,E,P^c\times_\delta\C\textmap{\simeq}\det E)$, and let
$\Sigma^{\pm}$  be the spinor bundles associated with $P^c$.  The spinor bundles
$\Sigma^{\pm}(P^h)$ of
$P^h$ and the automorphism-bundle $\G(P^h)$ can be be expressed
in terms of the
$G_3$-reduction as follows
$$\Sigma^{\pm}(P^h)=[\Sigma^{\pm}]^{\vee}\otimes E=
\Sigma^{\pm}\otimes E^{\vee}\ ,\ \
\G(P^h)=SU(E) \ .
$$

Moreover, the associated $PU(2)$-bundle $\delta(P^h)=
P^h\times_\delta PU(2)$ is naturally
isomorphic to the $S^1$-quotient of the unitary frame bundle $P_E$ of $E$.
\vspace{0.5cm}\\ \\
3.   $Spin^{U(2)}$-structures: \\

Consider the $U(2)$ spin   group
$$Spin^{U(2)}:=Spin\times_{\Z_2} U(2)\ ,$$
and let $p:U(2)\map  PU(2)$ be the canonical projection. The map
$$p\times\det: U(2)\map PU(2)\times S^1$$
 induces an isomorphism
$\qmod{U(2)}{\{\pm\id\}}=PU(2)\times S^1$. Therefore the map
$\delta:Spin^{U(2)}\map
\qmod{U(2)}{\{\pm\id\}}$   can be written  as a pair $(\bar\delta,\det)$
consisting of a
$PU(2)-$ and an
$S^1$-valued morphism. We have  exact sequences
\begin{equation}
\begin{array}{c}1\map Spin \map Spin^{U(2)}\textmap{(\bar\delta,\det)}
PU(2)\times
S^1\map  1 \\
\\
1\map  U(2)\map  Spin^{U(2)}\textmap {\pi} SO \map  1 \\ \\
 1\map \Z_2\map Spin^{U(2)} \textmap{(\pi,\bar\delta,\det)} SO
\times PU(2)\times S^1
\map 1 \\ \\
1 \map SU(2)\map Spin^{U(2)}\textmap{(\pi,\det)} SO \times S^1 \map 1 \ .
\end{array}
\end{equation}

Let $P^u\map P $ be a  $Spin^{U(2)}$-structure  in a $SO$-bundle $P$  over $X$.

An important role will be played by the subbundles
$$\G_0(P^u):=P^u\times_{\Ad_{U(2)}}
SU(2)\ ,\ \ \gr_0(P^u):=P^u\times_{\Ad_{U(2)}}su(2)$$
of $\G(P^u)=P^u\times_{\Ad_{U(2)}}U(2)$,
$\gr(P^u):=P^u\times_{\Ad_{U(2)}} u(2)$ respectively. The group of sections
$${\cal G}_0(P^u):=\Gamma(X,\G_0(P^u))$$
in   $\G_0(P^u)$ can be identified with the group of automorphisms of $P^u$ over
the $SO\times S^1$-bundle   $P\times_X (P^u\times_{\det} S^1)$.

By Propositions 1.1.4, 1.1.5 we get
\begin{pr} Let $P$ be a principal $SO$-bundle, $\bar P$ a $PU(2)$-bundle, and
$L$ a Hermitian line bundle  over $X$.\\
 i) $P$ admits a $Spin^{U(2)}$-structure $P^u
\rightarrow  P$  with
$$P^u\times_{\bar \delta}PU(2)\simeq\bar P\ ,\ \ P^u\times_{\det}\C\simeq L$$
iff $w_2(P)=w_2(\bar P)+\overline c_1(L)$, where $\overline c_1(L)$ is
the mod 2 reduction of $c_1(L)$ .\\
 ii) If the base $X$ is a compact oriented 4-manifold, then the
map
$$P^u\longmapsto \left([P^u\times_{\bar\delta} PU(2)]
,[P^u\times_{\det}\C]\right)$$
 defines a 1-1 correspondence between the set of
isomorphism classes of
 $Spin^{U(2)}$-struc\-tures in  $P$  and   the set of pairs of isomorphism
 classes
$([\bar P],[L])$, where
$\bar P$ is a $PU(2)$-bundle and $L$ an $S^1$-bundle with $w_2(
P)=w_2(\bar P)+\overline c_1(L)$. The latter set can be identified with
$$\{(p,c)\in H^4(X,\Z)\times H^2(X,\Z) |\ p\equiv (w_2(P)+ \bar c)^2\ {\rm
mod}\ 4\}
$$
\end{pr}
\qed

The group $Spin^{U(2)}(4)=\qmod{SU(2) \times SU(2)
\times U(2)}{\{\pm(\id,\id,\id)\}}$
fits in the exact sequence
$$1\map  S^1 \map \tilde G_3\map  Spin^{U(2)}(4)\map  1\ ,$$
where
$$\tilde G_3:=\{(a,b,c)\in U(2)\times U(2)\times U(2)|\ \det a=\det b\}\ .$$
and a $Spin^{U(2)}(4)$-bundle $P^u$ admits a $\tilde G_3$-reduction iff
$w_2(P^u\times_\pi SO(4))$ has integral lifts. Therefore, as in
 Proposition 1.1.9, we get

\begin{pr} Let $P$ be an $SO(4)$-bundle whose second Stiefel-Whitney
 class admits
integral lifts.

There is a 1-1 correspondence between isomorphism classes of
  $Spin^{U(2)}$-structures   in $P$ and   equivalence classes of pairs
consisting of a $Spin^c(4)$-structure $P^c\map  P$ in
$P$ and a  $U(2)$-bundle $E$. Two pairs are considered equivalent if,
 after tensoring
the first one with a line bundle, they become isomorphic over $P$.
\end{pr}

Suppose that the $Spin^{U(2)}(4)$-bundle $P^u$ admits an $\tilde
G_3$-lifting,    let
$(P^c,E)$ be the pair associated with this reduction, and let $\Sigma^{\pm}$
 be the spinor
bundles associated with $P^c$. Then  the  associated bundles
$\Sigma^{\pm}(P^u)$,
$\bar\delta(P^u)$,
$\det(P^u)$, ${\G(P^u)}$, $\G_0(P^u)$
  can be expressed in terms of the pair $(P^c,E)$ as follows:
$$\Sigma^{\pm}(P^u)=[\Sigma^{\pm}]^{\vee}\otimes
E=\Sigma^{\pm}\otimes (E^{\vee}\otimes[\det(P^u)]) \ ,\ \
\bar\delta(P^u)\simeq
\qmod{P_E}{S^1}\ , $$ $$  \det(P^u)\simeq \det (P^c)^{-1}\otimes (\det E)\ , \ \
\G(P^u)=U(E),\ \gr(P^u)=u(E)\ ,$$
$$ \ \G_0(P^u)=SU(E),\ \gr_0(P^u)=su(E)\ .
 $$

  \subsection{The G-monopole equations}
\subsubsection{Moment maps for families of complex structures}

Let $(M,g)$ be a   Riemannian manifold, and ${\cal J}\subset
A^0(so(T_M))$ a family of complex structures on
$M$ with the property that $(M,g,J)$ is a K\"ahler manifold, for
 every $J\in {\cal J}$. We
denote by $\omega_J$ the K\"ahler form of this K\"ahler manifold.
 Let $G$ be a compact
Lie group acting on $M$ by isometries with are holomorphic with
respect to any $J\in{\cal
J}$. Let $U$ be a fixed   subspace of  $A^0(so(T_M))$ containing the
family ${\cal J}$, and suppose for simplicity that $U$ is finite dimensional. We
define the \underbar{total} \underbar{K\"ahler} \underbar{form}
$\omega_{\cal J}\in
A^2(U^{\vee})$ by the formula
$$\langle\omega_{\cal J},u\rangle=g(u(\cdot),\cdot)\ .$$
\begin{dt} Suppose  that the total K\"ahler form is closed and
$G$-invariant.
   A map $\mu:M \map \Hom(\g, U^{\vee})$ for will
 be called a ${\cal
J}$-moment map  for the
$G$-action in
$X$   if the following two identities hold \\
1. $\mu(ag)=(\ad_g\otimes\id_{U^{\vee}})(\mu(a))$ $\forall\ a\in M,\ g\in G$.\\
2. $d(\langle\mu,{\alpha}\rangle)=\iota_{\alpha^{\#}}\omega_{\cal J}$ in
$A^1(U^{\vee})$
$\forall\
\alpha\in \g$, where $\alpha^{\#}$ denotes the vector field associated
with $\alpha$.
\end{dt}

In many cases $\g$ comes with a natural $\ad$-invariant euclidean metric.
 A map $\mu:M
\map \g\otimes U$ will be   called also a moment map if its composition
with the
morphism  $\g\otimes U\map \g^{\vee}\otimes U^{\vee}$ defined by the
euclidean structures
in $\g$ and $U$ is a moment map in the above sense. Similarly, the total
K\"ahler form can be
regarded (at least in the finite dimensional case) as an element in
$\omega_{\cal J}\in
A^2 (U)$.

Note that if $\mu$ is a moment map with respect to ${\cal J}$, then
for every
$J\in{\cal J}$ the map $\mu_J:=\langle \mu,J\rangle:M\map\g^{\vee}$ is a
moment map for
the $G$-action in $X$ with respect to the symplectic structure
$\omega_J$.
\

\begin{re} Suppose
that the total K\"ahler form $\omega_{\cal J}$ is $G$-invariant and closed.
Let $\mu:M\map
\Hom(\g, U^{\vee})$ be a ${\cal J}$-moment map    for a
free $G$-action and suppose that $\mu$ is a submersion at every point in
$\mu^{-1}(0)$.
Then
$\omega_{\cal J}$ descend to a closed
$U^{\vee}$-valued 2-form  on the quotient
manifold $\qmod{\mu^{-1}(0)}{G}$. In particular,
in this case, all the 2-forms
$\omega_J$ descend to closed   2-formes on this quotient, but they may be
degenerate.

\end{re}
\vspace{3mm}
{\bf Examples:}\hfill{\break}\\
1. Hyperk\"ahler manifolds: \\

Let $(M,g,(J_1,J_2,J_3))$ be a hyperk\"ahler manifold [HKLR]. The
three complex structures
$J_1,J_2,J_3$ span a sub-Lie algebra $U\subset A^0(so(T_M))$
naturally isomorphic to
$su(2)$. Suppose for simplicity that $Vol(M)=1$. The sphere $S(U,\sqrt
2)\subset U$ of
radius
$\sqrt 2$ contains  the three complex structures and for any
$J\in  S(U,\sqrt 2)$ we get a K\"ahler manifold $(M,g,J)$. Suppose  that
$G$ acts on $M$ preserving the hyperk\"ahler structure. A hyperk\"ahler
moment map $\mu:M\map
\g\otimes su(2)$ in the sense  of   [HKLR] can be regarded as a moment
map with respect to the family  $S(U,\sqrt 2)$ in the sense above. If the
assumptions in the Remark above are fulfilled, then the forms
$(\omega_J)_{J\in S(U,\sqrt 2)}$ descend to
\underbar{symplectic} forms on the quotient  $\qmod{\mu^{-1}(0)}{G}$,
which can be endowed with a natural hyperk\"ahler structure in this
way [HKLR].\\ \\
%
2.  Linear  hyperk\"ahler spaces:\\

Let $G$ be a compact Lie group  and
$G\subset U(W)$   a unitary representation of
$G$.  A moment map for the
$G$-action on $W$ is given by
$$\mu_G(w)=-\Pr_\g\left(\frac{i}{2}(w\otimes\bar w)\right)
$$
where $\Pr_\g:u(W)\map \g=\g^{\vee}$ is the projection $\g\hookrightarrow
u(V)$. Any other moment map can be obtained by adding a constant central
element in $\g$.

In the special case of   the standard left action of $SU(2)$ in $\C^2$, we
denote by $\mu_0$
the associated moment map. This is given by

$$\mu_0(x)=-\frac{i}{2}(x\otimes\bar x)_0 \ ,
$$
where  $(x\otimes\bar x)_0$ denotes the trace-free component of the Hermitian
endomorphism $x\otimes\bar x$.

Consider now the scalar  extension
$M:=\H\otimes_{\C} W$. Left multiplications by  quaternionic units define a
$G$-invariant
hyperk\"ahler structure in $M$. The corresponding family of complex
structures is
parametrized by the radius $\sqrt 2$-sphere $S$ in the space of imaginary
quaternions
identified with $su(2)$.

Define the quadratic map $\mu_{0,G}:\H\otimes_\C W\map su(2)\otimes \g$ by
$$\mu_{0,G}(\Psi)=\Pr_{[su(2)\otimes \g]}(\Psi\otimes\bar\Psi) \ .
$$
It acts on tensor monomials by
$$x\otimes w\stackrel{\mu_{0G}}{\longmapsto} -4\mu_0(x)\otimes\mu_G(w)\in
su(2)\otimes \g\subset \Herm(\H\otimes_\C W)\ .
$$

It is easy to see that $-\frac{1}{2}\mu_{0,G}$ is a  moment map for the
$G$ action in $M$ with respect to the linear hyperk\"ahler structure in $M$
introduced above.
\\ \\
3.  Spaces of spinors:\\

Let $P^G$ be $Spin^G(4)$-bundle over a compact Riemannian manifold $(X,g)$.
The corresponding spinor bundles $\Sigma^{\pm}(P^G)$ have
$\H_{\pm}\otimes_\C V$ as  standard fibres. Any section $J\in
\Gamma(X,S(\ad_{\pm}(P^G),\sqrt 2))$ in the radius $\sqrt 2$-sphere bundle
associated to
  $ad_{\pm}(P^G)$ gives a complex (and hence a K\"ahler) structure in
$A^0(\Sigma^{\pm}(P^G))$.

Therefore (after suitable Sobolev completions)
the space of
sections  $$\Gamma(X,S(\ad_{\pm}(P^G),\sqrt 2))$$ can be regarded as a family
of K\"ahler
structures in  the  space of sections $A^0(\Sigma^{\pm}(P^G))$ endowed with
the standard
$L^2$-Euclidean metric.  Define a quadratic map $\mu_{0,{\cal
G}}:A^0(\Sigma^{\pm}(P^G))\map A^0(ad_{\pm}(P^G)\otimes \gr)$ by sending an
element
$\Psi\in A^0(\Sigma^{\pm}(P^G))$ to the section in
$ad_{\pm}(P^G)\otimes\gr$ given by
the fibrewise projection of  $\Psi\otimes\bar\Psi\in
A^0(\Herm(\Sigma^{\pm}(P^G)))$.

Then $-\frac{1}{2}\mu_{0,{\cal G}}:A^0(\Sigma^{\pm}(P^G))\map
A^0(\ad_{\pm}(P^G)\otimes
\gr)\subset\Hom(A^0(\gr), A^0(\ad_{\pm}(P^G)^{\vee})$ can be
regarded as a   $\Gamma(X,S(\ad_{\pm}(P^G),\sqrt 2))$-moment  map for the
natural action of the
gauge group ${\cal G}$.   \\
\\
4. Spaces of connections on a 4-manifold:\\

Let $(X,g)$ be a compact oriented Riemannian 4-manifold, $G\subset U(r)$ a
compact Lie
group, and $P$ a principal $G$-bundle over $X$. The space of connections
${\cal A}(P)$ is an
euclidean affine space modelled on $A^1(\ad(P))$, and the gauge group ${\cal
G}:=\Gamma(X,P\times_{Ad}G)$ acts from the left by
$L^2$-isometries.  The space of
almost complex structures in $X$ compatible with the metric and the
orientation can be
identified with space of sections in the sphere bundle
$S(\Lambda^2_+,\sqrt 2)$ under the map which associates to an almost
complex structure
$J$ the K\"ahler form
$\omega_J:=g(\cdot,J(\cdot))$ [AHS]. On the other hand any almost complex
structure
$J\in \Gamma(X,S(\Lambda^2_+,\sqrt 2))$ induces a   gauge invariant
\underbar{integrable}  complex structure in the affine space ${\cal A}(P)$
by identifying
$A^1(\ad(P))$ with $A^{01}_J(\ad(P)^{\C})$.

The total K\"ahler form of this family is the element
$\Omega\in A^2_{{\cal A}(P)}(A^2_{+,X})$    given by
$$\Omega(\alpha,\beta) = \tr(\alpha\wedge\beta)^+ \ ,
$$
where $\alpha$, $\beta\in A^1(\ad(P))$.

Consider the map $  F^+ :{\cal A}(P)\map
A^0(\ad(P)\otimes \Lambda^2_{X,+})\subset  \Hom(A^0(\ad(P)),(A^2_+)^{\vee})$
given by
$A\longmapsto   F_A^+$.
It satisfies the equivariance property 1. in Definition 1.2.1. Moreover,
for every
$A\in{\cal A}(P)$,
$\alpha\in A^1(\ad(P))=T_A({\cal A}(P))$,
$\varphi\in A^0(\ad(P))=Lie({\cal G})$ and $\omega\in A^2_+$  we have
(denoting by $\delta$
the exterior derivative on ${\cal A}(P)$)
$$\left\langle(\iota_{\varphi^{\#}}\Omega)(\alpha)-
\langle\delta_A(F^+)
(\alpha),\varphi\rangle,\omega\right\rangle
=\langle
d^+[\tr(\varphi\wedge\alpha)],\omega\rangle=
\int_X\tr(\varphi\wedge\alpha)\wedge
d\omega \ .
$$
This formula means that the second condition in Definition 1.2.1.  holds
up to 1-forms on ${\cal A}(P)$ with values in the subspace $\im[d^+: A^1_X
\rightarrow
A^2_{X,+} ]$. Let $\bar\Omega$ be the image of $\Omega$ in $A^2_{{\cal
A}(P)}\left[\qmod{A^2_{X,+}}{\im(d^+)}\right]$. Putting
$${\cal A}^{ASD}_{reg}=\{A\in{\cal A}(P)|\ F_A^+=0,\ {\cal G}_A=Z(G),\
H^0_A=H^2_A=0\}$$
we see that   $\bar\Omega$   descends to a closed
$\left[\qmod{A^2_{X,+}}{\im(d^+)}\right]$-valued 2-form
$[\bar\Omega]$ on the moduli space of regular anti-selfdual connections
${\cal M}^{ASD}_{reg}:=\qmod{{\cal A}^{ASD}_{reg}}{{\cal G}}$.  Thus we
may consider  the
map
 $F^+$  as a   $\Gamma(X,S(\Lambda^2_+,\sqrt 2))$-moment map modulo $d^+$-exact
forms for the  action of the gauge group on ${\cal A}(P) $.

Note that in the case $G=SU(2)$ taking $L^2$-scalar product of
$\frac{1}{8\pi^2}[\bar\Omega]$ with  a harmonic selfdual form $\omega\in
\H^2_+$   defines a de Rham representant of Donaldson's
$\mu$-class associated with the Poncar\'e dual of $[\omega]$.

The following simple consequence of the above observations can be regarded
as the
starting point of Seiberg-Witten theory.

\begin{re} The data of a $Spin^G(4)$-structure in the Riemannian manifold
$(X,g)$ gives an
isometric  isomorphism $\frac{1}{2}\Gamma:\Lambda^2_{+}\map \ad_{+}(P^G)$. In
particular we get an identification between the two familes
$\Gamma(X,S(\ad_{+}(P^G),\sqrt 2))$ and $\Gamma(X,S(\Lambda^2_{+},\sqrt 2))$ of
complex structures in $A^0(\Sigma^{+}(P^G))$ and ${\cal A}(\delta(P^G))$
studied
before.  Consider the action of the gauge group  ${\cal G}:=\Gamma(X,\G)$ on
the product ${\cal A}(\delta(P^G))\times A^0(\Sigma^{+}(P^G))$ given by
$$[(A,\Psi),f]\longmapsto (\delta(f)(A),f (\Psi))\ .
$$

 This action admits a (generalized)
moment map modulo
$d^+$-exact  forms (with respect to the family
$\Gamma(X,S(\ad_{+}(P^G),\sqrt 2))$)  which
is given by the formula
$$(A,\Psi)\longmapsto   F_A^{+}-\Gamma^{-1}(\mu_{0,{\cal G}}(\Psi)) \ .
$$
\end{re}

\subsubsection{Dirac harmonicity and the $G$-monopole equations}

Let $P^G$ be a $Spin^G$-bundle. Using the third exact sequence in $(*)$
 sect. 1.1, we see
that the data of a connection in $P^G$ is equivalent to the data of a pair
consisting of a
connection in
$P^G\times_\pi SO$, and a connection in $\delta(P^G)$. In particular, if
$P^G\map  P_g$ is
 a $Spin^G(n)$-structure in the frame bundle of an oriented Riemannian
$n$-manifold $X$,
then the data of a connection $A$ in    $\delta(G)$  is equivalent to the
data of a
connection
$B_A$ in $P^G$ lifting the Levi-Civita connection in $P_g$. Suppose now
 that $n=4$, and
denote as usual  by $\gamma: \Lambda^1\map \Hom_\G(\Sigma^+(P^G)^+,
\Sigma^-(P^G))$
 the Clifford map of   a fixed
$Spin^G(4)$-structure $P^G\stackrel{\sigma}\map P_g$, and by
$\Gamma:\Lambda^2_{\pm}\map
\ad_{\pm}(P^G)$ the induced isomorphisms. We define the Dirac
 operators
$\Dr_A^{\pm}$ associated with $A\in{\cal A}(\delta(P^G))$ as the
composition
$$A^0(\Sigma^{\pm}(P^G))\textmap{\nabla_{{B_A}}}
A^1(\Sigma^{\pm}(P^G))\textmap{\cdot\gamma} A^0(\Sigma^{\mp}(P^G)) \ .
$$
We put also
$$\Sigma(P^G):=\Sigma^+(P^G)\oplus \Sigma^-(P^G)\ ,\ \
\Dr_A:=\Dr_A^+\oplus\Dr_A^-:A^0(\Sigma(P^G))\map A^0(\Sigma(P^G))\ .$$
Note that $\Dr_A$ is a  self-adjoint first order elliptic operator.
\begin{dt} A pair $(A,\Psi)\in{\cal A}(P^G)\times A^0(\Sigma(P^G))$
will be called
(Dirac) harmonic if $\Dr_A\Psi=0$.
\end{dt}

The harmonicity condition is obviously   invariant with respect to the
 gauge group ${\cal
G}(P^G):=\Gamma(X,\G(P^G))$. The monopole equations associated to
 $\sigma$ couple the
two gauge  invariant equations we introduced above: the vanishing of
 the "moment map "
(cf. 1.4.1) of the gauge action with respect to the family of complex
 structures
$\Gamma(X,S(\ad_+(P^G),\sqrt 2))$ in the affine space
${\cal A}(\delta(P^G))\times A^0(\Sigma^+(P^G))$ and the Dirac
harmonicity.

\begin{dt} Let $P^G\textmap{\sigma}  P$ be a $Spin^G(4)$-structure on
the compact
oriented Riemannian 4-manifold $X$.
The  associated Seiberg-Witten equations for a pair
$(A,\Psi)\in {\cal A}(\delta(P^G)) \times A^0(\Sigma^+(P^G))$   are
$$\left\{\begin{array}{ccc}
\Dr_A\Psi&=&0\\
\Gamma(F_A^+)&=&\mu_{0,{\cal G}}(\Psi)
\end{array}\right. \eqno{(SW^\sigma)}$$
\end{dt}

The solutions of these equations modulo the gauge group will be
 called $G$-monopoles.

The case $G=S^1$ corresponds to the classical (abelian)
Seiberg-Witten theory. The case
$G=SU(2)$ was extensively studied in [OT5], and from a physical
 point view in [LM].

 \begin{re} If the Lie algebra $\g$ of $G$ has non-trivial center $z(\g)$,
then the moment
map of the  gauge action in   $A^0(\Sigma^+(P^G))$ is not unique. In this
case it is
more natural to consider the  family of equations
$$\left\{\begin{array}{ccc}
\Dr_A\Psi&=&0\\
\Gamma(F_A^+)&=&\mu_{0,{\cal G}}(\Psi)+\beta \ ,
\end{array}\right. \eqno{(SW^\sigma_\beta)}$$
 obtained by adding in the  second equation a
section
$$\beta\in A^0(\ad_+(P^G)\otimes z(\g))\simeq A^2_+(X,z(\g))\ .$$
\end{re}

In the case $G=S^1$ the  equations of this form are called \underbar{twisted}
\underbar{monopole} equations [OT6].   If
$b_+(X)=1$, the   invariants defined using moduli spaces of twisted
monopoles depend in an
essential way  on the twisting term $\beta$ ([LL], [OT6]).

The particular case $G=U(2)$ requires a separate discussion, since
in this case
$\delta(U(2))\simeq PU(2)\times S^1$ and, correspondingly, the bundle
 $\delta(P^u)$
associated with a $Spin^{U(2)}(4)$-structure $P^u\textmap{\sigma}  P_g$
splits as the
product
$$\delta(P^u)=\bar\delta(P^u)\times_X\det(P^u)$$
 of a $PU(2)$-bundle with a $U(1)$-bundle. The  data
of a connection in $P^u$ lifting the Levi-Civita connection in $P_g$ is
 equivalent to the data
of a pair
$A=(\bar A,a)$ formed by a connection $\bar A$ in $\bar\delta(P^u)$ and a
 connection $a$ in
$\det(P^u)$. An alternative approach regards the connection
$a\in {\cal A}(\det(P^u))$ as a
parameter   (not an unknown !) of the equations, and studies the
 corresponding monopole
equations for a pair $(\bar A,\Psi)\in {\cal A}(\bar\delta(P^u))\times
 A^0(\Sigma^+)$.
$$\left\{\begin{array}{ccc}
\Dr_{\bar A,a}\Psi&=&0\\
\Gamma(F_{\bar{A}}^+)&=&\mu_{0,0}(\Psi)
\end{array}\right. \eqno{(SW^\sigma_a)}$$
Here $\Dr_{\bar A,a}$ denotes the Dirac operator associated to the
 connection in $P^u$
which lifts the Levi-Civita connection in $P_g$, the connection $\bar A$
in the
$PU(2)$-bundle
$\bar\delta(P^u)$ and the connection $a$ in the $S^1$-bundle $\det P^u$;
the quadratic
map $\mu_{0,0}$ sends a spinor $\Psi\in A^0(\Sigma^+(P^u))$  to the
projection of the
endomorphism
$(\Psi\otimes\bar\Psi)\in A^0(\Herm(\Sigma^+(P^u)))$ on
$A^0(\ad_+(P^u)\otimes\gr_0(P^u))$.

  The natural gauge group
which lets invariant the equations is the group ${\cal G}_0(P^u):=
\Gamma(X,\G_0(P^u))$ of
automorphisms of the bundle $P^u$ over the bundle-product $P_g\times_X
\det(P^u)$,and
$-\mu_{0,0}$ is the
$\Gamma(X,S(\ad_+,\sqrt 2))$-moment map for the ${\cal G}_0(P^u)$-action in the
configuration space.  There is no ambiguity in choosing the moment map of the
${\cal G}_0(P^u)$-action, so there is \underbar{no} natural way to
perturb   these equations besides varying the connection-parameter $a\in{\cal
A}(\det(P^u))$.

 Since the connection-component   of the unknown is a
$PU(2)$-connection,  these  equations will   be called the
$PU(2)$-\underbar{monopole}      \underbar{equations}, and its solutions
 modulo the
gauge group ${\cal G}_0(P^u)$ will be called
$PU(2)$-monopoles.

Note that if the $Spin^{U(2)}(4)$-structure $P^u\map P_g$ is associated
with the pair $(P^c
\map P_g,E)$ (Proposition  1.1.11), the quadratic map $\mu_{0,0}$ sends
 a spinor $\Psi\in
A^0\left(\Sigma^+(P^c)\otimes [E^{\vee}\otimes\det (P^u)]\right)$  to the
 projection of
$$(\Psi\otimes\bar\Psi)\in
A^0\left(\Herm\left(\Sigma^+(P^c)
\otimes[E^{\vee}\otimes\det (P^u)]\right)\right)$$
 on
$A^0\left(su(\Sigma^+)\otimes su([E^{\vee}
\otimes\det (P^u)]\right)$.

\begin{re} The data of a $Spin^h(4)$-structure in $X$ is equivalent to
 the data of
$Spin^{U(2)}$-structure $P^u\textmap{\sigma} P$ together with a
trivialization of the
$S^1$-bundle $\det(P^u)$. The corresponding $SU(2)$-Seiberg-Witten
 equations coincide with
the $PU(2)$-equations $SW^\sigma_\theta$ associated with the trivial
 connection
$\theta$  in the trivial bundle  $\det(P^u)$.
\end{re}

We shall always regard the $SU(2)$-monopole equations as
special $PU(2)$-monopole
equations. In particular we shall use the notation
$\mu_{0,{\cal G}}=\mu_{0,0}$ if ${\cal
G}$ is the gauge group associated with a $Spin^h(4)$-structure.

\begin{re} The moduli space of
$PU(2)$-monopoles of the form $(\bar A,0)$ can be identified
 with a moduli space of
anti-selfdual
$PU(2)$-connections, modulo the gauge group ${\cal G}_0$. The
 natural morphism of
${\cal G}_0$ into the usual $PU(2)$-gauge group of
bundle automorphisms of $\bar\delta(P^u)$ is a local isomorphism
 but in general it is
not surjective (see [LT] ). Therefore the space of $PU(2)$-monopoles
 of the form $(\bar
A,0)$ is a finite cover of the corresponding Donaldson moduli space of
 $PU(2)$-instantons.
\end{re}

\begin{re} Let $G$ be a compact Lie group endowed with a central invlotion
$\iota$ and an
arbitrary unitary representation $\rho:G\map  U(V)$ with
$\rho(\iota)=-\id_V$. One can
associate to any $Spin^G(4)$-bundle the spinor bundles $\Sigma^{\pm}$ of
standard fibre
$\H_{\pm}\otimes V$. Endow the Lie algebra
$\g$ with an
$\ad$-invariant metric. Then one can define
$\mu_{0,G}$ using the adjoint of the map $\g\map  u(V)$ instead of the
orthogonal projection, and
the
$G$-monopole equations have sense in this more general framework.
\end{re}

\subsubsection{Reductions}

Let $H\subset G\subset U(V)$ be a closed subgroup of $G$ with $-\id_V\in H$.
Let $P^G\textmap{\sigma}  P$ be a $Spin^G $-structure in the $SO$-bundle
bundle $P$.
\begin{dt} A $Spin^H$-reduction of $\sigma$ is a subbundle  $P^H$ of
$P^G$ with structure group $Spin^H\subset Spin^G$.
\end{dt}

Note that such a reduction $P^H\hookrightarrow P^G$ defines a reduction
$\delta(P^H)\hookrightarrow\delta(P^G)$ of the structure group of the
bundle $\delta(P^G)$
from $\qmod{G}{\Z_2}$ to
$\qmod{H}{\Z_2}$, hence it defines in particular an injective linear
morphism  ${\cal
A}(\delta(P^H))\hookrightarrow{\cal A}(\delta(P^G))$ between the associated
affine spaces
of connections.

Let now $V_0$ be an  $H$-invariant subspace of $V$.

Consider a $Spin^G(4)$-structure $P^G\textmap{\sigma} P$  in the
$SO(4)$-bundle $P$,  and
a  $Spin^H(4)$-reduction $P^H\stackrel{\rho}{\hookrightarrow} P^G$  of
$\sigma$. Let
$\Sigma^{\pm}(P^H,V_0)$ be the spinor bundles associated with $P^H$ and the
$Spin^H(4)$-representation in
$\H^{\pm}\otimes_\C V_0$.

The inclusion $V_0\subset V$ induces bundle inclusions of the associated
spinor bundles
$\Sigma^{\pm}(P^H,V_0)\hookrightarrow \Sigma^{\pm}(P^G)$. Suppose now that
$P_g$ is the frame-bundle of a compact oriented Riemannian 4-manifold,   choose
 $A\in{\cal
A}(\delta(P^H))\subset {\cal A}(\delta(P^G))$, and let be $B_A\in{\cal
A}(\delta(P^G))$ be
the
induced connection.   Then  the spinor bundles $\Sigma^{\pm}(P^H,V_0)$ become
$B_A$-parallel subbundles of
$\Sigma^{\pm}(P^G)$, and the Dirac operator
$$\Dr_A:\Sigma(P^G)\map  \Sigma(P^G)$$
maps $\Sigma(P^H,V_0)$ into itself. Therefore the set of Dirac-harmonic pairs
associated with
$(\sigma\circ\rho,V_0)$ can be identified with  a subset of the set of
Dirac-harmonic
pairs associated with $(\sigma,V)$.

The group $G$ acts on the set
$$\{(H,V_0)|\ H\subset G\ {\rm closed\ subgroup},\ V_0\subset V\ {\rm is}\
H-{\rm
invariant}\}\ .
$$
of subpairs of $(G,V)$ by $[g,(H,V_0)]\longmapsto( Ad_g(H),g(V_0))$.
Moreover, for any
$Spin^H(4)$-reduction
$P^H\hookrightarrow P^G$ of
$\sigma$ and any element $g\in G$ we get a reduction
$P^{\Ad_g(H)}\hookrightarrow P^G$ of
$\sigma$ and subbundles $\Sigma^{\pm}(P^{\Ad_g(H)},g(V_0))$ of the spinor
bundles
$\Sigma^{\pm}(P^G)$.

\begin{dt} A subpair $(H,V_0)$ of $(G,V)$ with  $-\id_V\in H$  will be
called admissible
and
$\mu_G|_{V_0}$ takes values in $\hg$ or, equivalently, if $\langle ik(v),
v\rangle=0$ for all
$k\in\hg^{\bot_{\g}}$ and $v\in V_0$.
\end{dt}
Therefore, if $(H,V_0)$ is admissible, then $\mu_G|_{V_0}$  can be
identified with the
  moment map $\mu_H$ associated with the $H$-action in $V_0$ (with
respect to the
metric in \hg induced from  \g -- see Remark 1.2.9).   If generally
$E$ is a system of  equations on a configuration space ${\cal A}$ we denote
by ${\cal
A}^E$ the space  of solutions of this system, enowed with the induced topology.
\begin{pr} Let $(H,V_0)$ be an \ admissible \ subpair of \ $(G,V)$.  \ A \\
$Spin^H(4)$-reduction
$P^H\textmap{\rho} P^G$   of the $Spin^G(4)$-structure $P^G\textmap{\sigma}
P_g$
induces an inclusion
$$\left[{\cal A}(\delta(P^H))\times
A^0(\Sigma^+(P^H))\right]^{SW^{\sigma\circ\rho}}\subset
\left[{\cal A}(\delta(P^G))\times A^0(\Sigma^+(P^G))\right]^{SW^{\sigma}}$$
which is equvariant with respect to the actions of the two gauge groups.
\end{pr}
\begin{dt} Let $(H,V_0)$ be an admissible subpair. A solution $(A,\Psi)\in
\left[{\cal
A}(\delta(P^G))\times A^0(\Sigma^+(P^G))\right]^{SW^\sigma}$ will be called
\underbar{reducible} \underbar{of} \underbar{type} $(H,V_0)$, if it belongs
to the image of
such an inclusion, for a suitable  reduction $P^H\textmap{\rho} P^G$.
\end{dt}

If $(H,V_0)$ is
admissible, $H\subset H'$ and $V_0$ is $H'$-invariant, then $(H',V_0)$ is
also admissible.
An admissible pair
$(H,V_0)$
will be called \underbar{minimal} if $H$ is  minimal in the set of closed
subgroups
$H'\subset G$ such that $(H',V_0)$ is an admissible subpair of $(G,V)$. The
sets of (minimal)
admissible pairs is
 invariant
under the natural $G$-action. We list the conjugacy classes of proper
minimal admissible
subpairs  in the  cases
$(SU(2),\C^{\oplus 2})=(Sp(1),\H)$,
$(U(2),\C^{\oplus 2})$,
$(Sp(2),\H^{\oplus 2})$.    Fix first   the maximal tori
$$
T_{SU(2)}:=\left\{\left(\matrix{z&0\cr 0& z^{-1}}\right)|z\in S^1\right\}\ ,\ \
T_{U(2)}:=\left\{\left(\matrix{u&0\cr 0& v }\right)| u,v\in S^1\right\}
$$
$$T_{Sp(2)}:=\left\{\left(\matrix{u&0\cr 0& v }\right)|u,v\in S^1\right\}
 $$
\\
On the right we list the minimal admissible subpairs  of the pair on the left:

$$\begin{array}{llcrl}
(SU(2),\C^{\oplus 2}):\ \ \ \ \ \ \ \ \ \ \ \ \ \ \ \ \ \ \   \  \ \ \   \
\ \ \ \ \ \ \ \ \ \
\ \ \ \ \   &(\{\pm1\}&, &\{0\}) &\\ \\ &(T_{SU(2)}&,&\C\oplus\{0\} )\\ \\
\end{array}
$$
$$\begin{array}{llcrl}
(U(2),\C^{\oplus 2}):\ \ \ \ \ \ \ \ \ \ \ \ \ \  \ \ \ \    \ \ \
&(\{\pm1\}&, &\{0\}) &\\ \\
&\left(\left\{\left(\matrix{\zeta&0\cr0&\pm1}\right)|\zeta\in
S^1\right\}\right.&,&\left.
\phantom{\matrix{1\cr1}}\C\times\{0\}\right) \\
\\
\end{array}
$$
$$\begin{array}{llcrl}
(Sp(2),\H^{\oplus 2}):\ \ \ \ \ \ \ \ \ &(\{\pm 1\}&, &\{0\}) &\\ \\
&\left(\left\{\left(\matrix{\zeta&0\cr0&\pm1}\right)|\zeta\in
T_{Sp(1)}\right\}\right.&,&\left.
\phantom{\matrix{1\cr1}}\C\oplus\{0_\H\}\right)\\
&\left(\left\{\left(\matrix{\zeta&0\cr0&\pm1}\right)|\zeta\in\ \  Sp(1)\
\ \right\}\right.&,&\left.
\phantom{\matrix{1\cr1}}\H\oplus\{0_\H\}\right)\\ \\
\end{array}
$$

\begin{re} Fix a maximal torus $T$ of $G$ with Lie algebra $\tg$,  and let
${\germ
W}\subset \tg^{\vee}$ be the weights of the induced $T$-action in $V$. Let
$V=\bigoplus\limits_{\alpha\in{\germ W}} V_\alpha$ be the corresponding
decomposition of
$V$ in  weight spaces.  If $(T,V')$ is a subpair of $(G,V)$, then $V'$ must
be a sum of
weight subspaces, i.e. there exist ${\germ W}'\subset {\germ W}$ such that
$V'=\bigoplus\limits_{\alpha\in{\germ W}'} V'_\alpha$, with $0\ne
V'_\alpha\subset
V_\alpha$.  When $G$ is one of the classical
groups  $SU(n)$, $U(n)$, $Sp(n)$ and $V$ the corresponding canonical
$G$-module, it follows
easily that $(T,V')$ is admissible iff $|{\germ W}'|=1$. Notice that there
is a natural action of
the Weil group $\qmod{N(T)}{T}$  in the set of abelian subpairs of the form
$(T,V')$.
\end{re}

 The case of the  $PU(2)$-monopole equations needs a separate discussion: Fix a
$Spin^{U(2)}(4)$-structure $\sigma:P^u\map  P_g$ in $P_g$ and a connection
$a$ in
 the
line bundle $\det (P^u)$.
In this case the admissible pairs are by definition equivalent to one of
$$(H,\{0\})\ ,\ \ H\subset U(2)\ {\rm with} -\id_V\in H\ ;\ \
(T_{U(2)},\C\oplus\{0\})
$$

An abelian reduction $P^{ T_{U(2)} }\stackrel{\rho}{\hookrightarrow} P^u$ of
$\sigma$  gives   rise to a pair of $Spin^c$-structures $(\cg_1:P^c_1\map  P_g,
\cg_2:P^c_2\map  P_g)$ whose determinant line bundles come with an
isomorphism
$\det(P^c_1)\otimes\det(P^c_2)=[\det (P^u)]^2$.  Moreover, the $PU(2)$-bundle
$\bar\delta(P^u)$ comes with an $S^1$-reduction $\bar\delta(P^u)=
P^{S^1}\times_\alpha
PU(2)$ where $[P^{S^1}]^2= \det(P^c_1)\otimes\det(P^c_2)^{-1}$ and
$\alpha:S^1\map
PU(2)$ is the standard embedding $\zeta\longmapsto\left[\left(\matrix{\zeta&0\cr
0&1}\right)\right]$. Since we have fixed the connection $a$ in  $\det (P^u)$,
the data of a
connection $\bar A\in{\cal A}(\bar\delta(P^u))$ which reduces to
$P^{ T_{U(2)} }$ via $\rho$ is equivalent to the data of a connection
$a_1\in{\cal A}(\det(P^c_1))$.

Moreover, we have  a natural parallel inclusion $\Sigma^{\pm}(P^c_1)\subset
\Sigma^{\pm}(P^u)$. Consider   the following twisted abelian monopole
 equations
[OT6] for  a pair $(A_1,\Psi_1)\in {\cal A}(\det (P^c_1))\times
A^0(\Sigma^{\pm}(P^c_1))$
$$\left\{\begin{array}{ccc}
\Dr_{A_1}\Psi_1&=&0\\
\Gamma(F_{A_1}^+)&=&(\Psi_1\bar\Psi_1)_0+\Gamma(F_a^+) \   .
\end{array}\right.  \eqno{(SW^{\cg_1}_{\Gamma(F_a^+)})}$$

Taking in Remark 1.2.6 as twisting term the form $\beta=\Gamma(F_a^+)$, we get
\begin{pr}A $Spin^{T_{U(2)}}$-reduction
$P^{ T_{U(2)}}\stackrel{\rho}{\hookrightarrow} P^u$
   of the $Spin^{U(2)}(4)$-structure $P^u\textmap{\sigma} P_g$
induces an inclusion
$$\left[{\cal A}(\det(P^c_1))\times
A^0(\Sigma^+(P^c_1))\right]^{SW^{\cg_1}_{\Gamma(F_a^+)}}\subset
\left[{\cal A}(\bar\delta(P^u))\times
 A^0(\Sigma^+(P^u))\right]^{SW^{\sigma}_a}$$
which is equivariant with respect to the actions of the two gauge groups.
\end{pr}

The fact that the  Donaldson ($PU(2)$-) $SU(2)$-moduli space is contained in the
 space of
($PU(2)$-) $SU(2)$-monopoles, and that
 (twisted) abelian monopoles arise as
abelian reductions in the space of ($PU(2)$)
$SU(2)$-monopoles suggests that these moduli spaces can be used to prove the
equivalence between the two theories [OT5].

This idea can be applied to get information about the Donaldson invariants
associated with
larger symmetry groups $G$ by relating these invariants to Seiberg-Witten
 type
invariants associated with smaller symmetry groups. In order to do this, one
 has first to
study invariants associated to the moduli spaces of reducible solutions of all
possible
types in a suitable moduli space  of $G$-monopoles.

\subsubsection{Moduli spaces of $G$-monopoles}

Let ${\cal A}$ be the configuration space
of one of the monopole equations $SW$ introduced in sect. 1.2.2.: For the
equations
$SW^\sigma_\beta$   associated with a $Spin^G(4)$-structure
$\sigma:P^G \map P_g$ in
$(X,g)$ and a section $\beta\in A^0(\ad_+(P^G)\otimes z(\g))$, the   space
${\cal A}$ coincides with  ${\cal A}(\delta(P^G))\times
A^0(\Sigma^+(P^G))$; in the
case of $PU(2)$-monopole equations $SW^\sigma_a$ associated to a
$Spin^{U(2)}(4)$-structure  $\sigma:P^u\map P_g$ and an abelian connection
$a\in{\cal
A}(\det(P^u))$ the configuration space is ${\cal A}(\bar\delta(P^u))\otimes
A^0(\Sigma^+(P^u))$. In this section, we denote by ${\cal G}$ the gauge group
corresponding
to the monopole equation $SW$, i.e. ${\cal G}={\cal G}(P^G)$ if
$SW=SW^\sigma_\beta$  and   ${\cal G}={\cal G}_0(P^u)$ in the $PU(2)$-case $SW=
SW^\sigma_a$. The Lie algebra $Lie({\cal G})$ of ${\cal G}$ is
$\Gamma(X,\gr(P^G))$
 in the
first case and $\Gamma(X,\gr_0(P^u))$  in the second.

The corresponding moduli space of $G$-monopoles  is defined as a
topological space by
$${\cal M}:=\qmod{{\cal A}^{SW}}{{\cal G}}\ .
$$

There is a standard way of describing the local structure of ${\cal M}$,
which was
extensively described in the cases $G=S^1$, $G=U(2)$ in [OT1] and in the
case $G=SU(2)$
(which is similar to the $PU(2)$-case) in [OT5] (see   [DK], [K], [LT], [M]
for the
instanton case and for the classical case of  holomorphic bundles).
We explain briefly the general strategy:

Let $p=(A,\Psi)\in {\cal A}^{SW}$. The infinitesimal action of
$Lie({\cal G})$
and  the
differential of $SW$ in $p$ define a "elliptic deformation complex"
$$0\map C^0_p \textmap{D^0_p} C^1_p\textmap{D^1_p}C^2_p \map 0
\eqno{({\cal C}_p)}
$$
where:\\
$C^0_p= Lie({\cal G})=\Gamma(X,\g(P^G))$  ( or $\Gamma(X,\g_0(P^u))$ in the
$PU(2)$-case),\\ \\
${\cal C}^1_p= T_p({\cal A})=A^1(\gr(P^G))\oplus A^0(\Sigma^{+}(P^G))$ (or
$A^1(\gr_0(P^u))\oplus A^0(\Sigma^{+}(P^u))$ in the
$PU(2)$-case),\\ \\
$C^2_p=A^0(\ad_+(P^G)\otimes\gr(P^G))\oplus A^0(\Sigma^{-}(P^G))$ (or
$A^0(\ad_+(P^u)\otimes\gr_0(P^u))\oplus A^0(\Sigma^{-}(P^u))$ in the
 $PU(2)$-case),\\ \\
 $D_p^0(f):=f^{\#}_p=(-d_A f, f\Psi)$,  \\ \\
$D_p^1(\alpha,\psi):=d_pSW(\alpha,\psi)=\left(\Gamma(d_A^+\alpha)-m
 (\psi,\Psi)-m
(\Psi,\psi),\gamma(\alpha)\Psi+\Dr_A(\psi)\right)\ .$
 Here $m$ is the
sesquilinear   map associated with the quadratic map $\mu_{0,{\cal G}}$ (or
$\mu_{0,0}$ in the
$PU(2)$-case).

The index $\chi$ of this elliptic complex is called the \underbar{expected}
\underbar{dimension} of the moduli space and can be easily computed by
Atiyah-Singer
Index-Theorem [LMi] in terms of   characteristic classes of
$X$ and vector bundles associated with $P^G$.

We give the result in the case of the $PU(2)$-monopole equations:
$$\chi(SW^\sigma_a)=\frac{1}{2}\left(-3 p_1(\bar\delta(P^u))+
c_1(\det(P^u))^2\right)-
\frac{1}{2}(3e(X)+4\sigma(X))
$$

The same methods as in [OT5] give:
\begin{pr}\hfill{\break}
1. The stabilizer ${\cal G}_p$ of $p$ is a finite dimensional Lie group
isomorphic to a
subgroup of $G$ which acts in a natural way in the harmonic spaces
$\H^i({\cal C}_p)$,
$i=0,\ 1,\ 2$.\\
2. There exists a neighbourhood $V_p$ of $P$ in ${\cal M}$, a ${\cal
G}_p$-invariant
neighbourhood $U_p$ of
$0$ in $\H^1({\cal C}_p)$, a  ${\cal G}_p$-equivariant real analytic map
$K_p:U_p\map
\H^2({\cal C}_p)$ with  $K_p(0)=0$, $dK_p(0)=0$ and a homeomorphism:
$$V_p\simeq \qmod{Z(K_p)}{{\cal G}_p}\ .
$$
\end{pr}

The homeomorphisms in the proposition above define a structure of a smooth
manifold of dimension $\chi$ in the open set
$${\cal M}_{reg}=\{[p]\in{\cal M}|\ {\cal G}_p=\{1\},\ H^2({\cal
C}_p)=0\}$$
 of regular points, and a structure of a real analytic orbifold in the open
set of points
with finite stabilizers.

Note that the stabilizer of a solution of the form  $(A,0)$ contains always
$\{\pm \id\}$,
hence   ${\cal M}$ has at least $\Z_2$-orbifold singularities in the
Donaldson points (see Remark 1.2.8).

As in the instanton case, the moduli space  ${\cal M}$ is in general
non-compact. The
construction of an Uhlenbeck-type compactification is treated in  [T1],
[T2].

\section{$PU(2)$-Monopoles and  stable oriented pairs}

In this section we show that the moduli spaces of $PU(2)$-monopoles  on a
compact
K\"ahler surface have  a natural complex geometric description in terms of
stable
oriented pairs.   We explain first briefly, following [OT5],
the concept  of oriented pair  and we indicate how  moduli space of   simple
oriented pairs are constructed. Next we restrict ourselves to the rank
2-case and we
introduce the concept of stable   oriented pair; the stability property we
need [OT5]
does \underbar{not} depend on a parameter and is an open  property.  An
algebraic
geometric approach can be found in [OST].

In section 2.2 we give a complex geometric description of the    moduli
spaces of
irreducible $PU(2)$-monopoles on  a K\"ahler surface in terms of  moduli
spaces of
stable oriented  pairs.  This description is used to give an explicit
description of
a moduli space of $PU(2)$-monopoles  on $\P^2$.

\subsection{Simple,  strongly simple and  stable oriented pairs}

Let $(X,g)$ be a compact K\"ahler manifold of dimension $n$, $E$ a
differentiable
vector bundle of rank
$r$ on
$X$, and ${\cal L}=(L,\bar\partial_{\cal L})$ a fixed holomorphic structure
in the
determinant line bundle $L:=\det E$.  We recall (see  [OT5]) the following
fundamental definition:
\begin{dt} An oriented pair of type $(E,{\cal L})$ is a pair $({\cal
E},\varphi)$,
where
${\cal E}$ is a holomorphic structure in $E$ such that $\det{\cal E}={\cal
L}$, and
$\varphi\in H^0({\cal E})$. Two oriented pairs $({\cal E}_1,\varphi_1)$, $({\cal
E}_2,\varphi_2)$ of type $(E,{\cal L})$ are called isomorphic if they are
congruent
modulo the natural action of the group $\Gamma(X,SL(E))$ of differentiable
automorphism of $E$ of determinant 1.
\end{dt}

Therefore we fix the underlying ${\cal C}^{\infty}$-bundle and the
holomorphic determinant line bundle (not only its isomorphism type !) of the
holomorphic bundles we consider.

An oriented pair  $p=({\cal E},\varphi)$  is called \underbar{simple} if its
 stabilizer
$\Gamma(X,SL(E))_p$ is contained in the center $\Z_r\id_E$ of
$\Gamma(X,SL(E))$, and is
called
\underbar{strongly} \underbar{simple} if its stabilizer is trivial.

The first property has an equivalent infinitesimal formulation: the pair
$({\cal
E},\varphi)$ is simple if and only if any trace-free holomorphic
endomorphism of
${\cal E}$ with
$f(\varphi)=0$ vanishes.

In [OT5] it was shown that
\begin{pr} There exists a  (possibly non-Hausdorff)  complex  analytic
orbifold   ${\cal M}^s(E,{\cal L} )$ parameterizing isomorphism classes of
simple oriented pairs of type
$(E,{\cal L})$. The open subset ${\cal M}^{ss}(E,{\cal L})\subset  {\cal
M}^{s}(E,{\cal L})$ consisting of strongly simple pairs is a complex
analytic space, and the points in ${\cal M}^s(E,{\cal L})\setminus{\cal
M}^{ss}(E,{\cal L})$ have neighbourhoods modeled on $\Z/r$-quotients.
\end{pr}

If  ${\cal E}$ is holomorphic bundle we denote by ${\cal S}({\cal E})$
the set of
reflexive subsheaves ${\cal F}\subset{\cal E}$ with $0<\rk({\cal
F})<\rk({\cal E})$.  Once we have fixed a section $\varphi\in H^0({\cal
E})$,      we
put
$${\cal S}_\varphi({\cal E}):=\{{\cal F}\in{\cal S}({\cal E})|\
\varphi\in H^0(X,{\cal F})\} \ .$$

We recall (see [B]) that ${\cal E}$ is called $\varphi$-\underbar{stable} if

$$\max (\mu_g({\cal E}),\sup\limits_{{\cal F}'\in{\cal S} ({\cal E})}
\mu_g({\cal F}'))<
\inf\limits_{{\cal F}\in {\cal S}_\varphi({\cal E})} \mu_g(\qmod{{\cal
E}}{{\cal F} })\ ,$$
where for a nontrivial torsion free coherent sheaf ${\cal F}$, $\mu_g({\cal
F})$
denotes its slope with respect to the K\"ahler metric $g$.  If the real number
$\lambda$ belongs to the interval $\left(\max (\mu_g({\cal E}),\sup
\limits_{{\cal
F}'\in{\cal S} ({\cal E})}
\mu_g({\cal F}')),
\inf\limits_{{\cal F}\in {\cal S}_\varphi({\cal E})} \mu_g(\qmod{{\cal
E}}{{\cal F} })\right)$, the pair $({\cal E},\varphi)$ is called
$\lambda$-stable.

If ${\cal
M}$ is a holomorphic line bundle and $\varphi\in H^0({\cal M})\setminus\{0\}$,
then $({\cal M},\varphi)$ is $\lambda$-stable iff $\mu_g({\cal M})<\lambda$.

The  correct definition of the stability property for oriented
pairs of arbitrary rank is a delicate point [OST]. The definition must
agree in the
algebraic-projective case with the corresponding GIT-stability condition.
On the
other hand, in the case $r=2$ the definition simplifies considerably and
this case is
completely sufficient for our purposes. Therefore from now on we assume $r=2$,
and we recall from [OT5] the following
\begin{dt} \hfill{\break}
An oriented pair $({\cal E},\varphi)$ of type $({ E},{\cal L})$
is called \underbar{stable} if one of the following conditions holds:\\
I.  \  ${\cal E}$ is
$\varphi$-stable, \\
II. $\varphi\ne 0$ and ${\cal E}$ splits in direct sum of line bundles
${\cal E}={\cal E}_1\oplus{\cal E}_2$, such that  \hspace*{5mm}  $\varphi\in
H^0({\cal E}_1)$ and  the pair $({\cal E}_1,\varphi)$ is $\mu_g({ E})$-stable.\\
   A holomorphic pair $({\cal E},\varphi)$ of type $({ E},{\cal L})$
is   called \underbar{polystable} if  it is stable, or $\varphi=0$  and
${\cal E}$
is a polystable bundle.
\end{dt}
\begin{re} An oriented  pair $({\cal E},\varphi)$ of type $(E,{\cal L})$ with
$\varphi\ne 0$ is stable iff     $\mu_g({\cal
O}_X(D_\varphi))<\mu_g(E)$, where $D_\varphi$ is the divisorial component of the
vanishing locus $Z(\varphi)$. An oriented pair of the form $({\cal E},0)$
is stable iff the
holomorphic  bundle ${\cal E}$ is stable.
\end{re}

\subsection{The projective vortex equation and stability of oriented pairs}

The stability property for holomorphic bundles has a well known  differential
geometric characterization: an holomorphic bundle is  stable if and only if
it is
simple and admits a Hermite-Einstein metric (see for instance [DK], [LT]).
Similarly,
an holomorphic pair $({\cal E},\varphi)$   is $\lambda$-stable if and only
it is simple
and ${\cal E}$ admits a Hermitian metric satisfying the vortex equation
associated
with the constant
$t=\frac{4\pi\lambda}{Vol_g(X)}$ [B]. All these important results are infinite
dimensional extensions of the {\it metric characterization of stability}
(see   [MFK],
[DK]).

The same approach gives in the case of oriented pairs the following
differential
geometric interpretation of stability [OT5]:

Let $E$ be a differentiable rank 2 vector bundle over a compact K\"ahler
manifold $(X,g)$, ${\cal L}$ a holomorphic structure in $L:=\det(E)$ and $l$ a
fixed Hermitian metric  in $L$.
\begin{thry} An holomorphic pair $({\cal E},\varphi)$   of type
$(E,{\cal L})$ with
$\rk({\cal E})=2$  is polystable iff ${\cal E}$ admits a
Hermitian metric $h$ with $\det h=l$ which solves  the following
\underbar{projective} \underbar{vortex} \underbar{equation}:
$$i\Lambda_g F_h^0 +\frac{1}{2}(\varphi\bar\varphi^h)_0=0\ .\eqno{(V)}$$
If $({\cal E},\varphi)$ is stable, then the
metric $h$ is unique.
\end{thry}

\begin{re} With an appropriate  definition of (poly)stability of oriented
pairs [OST], the
theorem holds for arbitrary rank $r$.
\end{re}

Denote by $\lambda\in{\cal A}(L)$ the the Chern connection of ${\cal L}$
associated with the metric $l$. Let $\bar{\cal A}_{\bar\partial_\lambda}$ be the
space of semiconnections in $E$ which induce the fixed semiconnection
${\bar\partial_\lambda}$ in $L$.

Fix a Hermitian metric $H$ in $E$ with $\det H=l$ and denote by   ${\cal
A}_\lambda$ the space of unitary connections in $E$ with induce the fixed
connection $\lambda$ in $L$. There is an obvious identification
${\cal A}_\lambda\textmap{\simeq}\bar{\cal A}_{\bar\partial_\lambda}$,
$C\longmapsto \bar\partial_C$ which endows the affine space  ${\cal A}_\lambda$
with a complex structure compatible with the standard $L^2$ euclidean structure.
Therefore, after suitable Sobolev completions,  the product ${\cal
A}_\lambda\times A^0(E)=\bar {\cal A}_{\bar\partial_\lambda}\times A^0(E)$
becomes a Hilbert K\"ahler manifold. Let ${\cal G}_0:=\Gamma(X,SU(E))$ be the
gauge group of unitary automorphisms of determinant 1 in
$(E,H)$ and let ${\cal G}_0^\C:=\Gamma(X,SL(E))$ be its complexification.
\begin{re} The map $m:{\cal A}_\lambda\times A^0(E)\map A^0(su(E))$ defined by
$$m(C,\varphi)=\Lambda_g F_C^0 -\frac{i}{2}(\varphi\bar\varphi^H)_0
$$
is a moment map for the ${\cal G}_0$-action in the K\"ahler manifold ${\cal
A}_\lambda\times A^0(E)$
\end{re}

If ${\cal E}$ is a holomorphic structure in $E$ with $\det{\cal E}={\cal L}$  we
denote by $C_{\cal E}\in {\cal A}_\lambda$ the Chern connection defined be
${\cal
E}$  and the fixed metric $H$.
The map $({\cal E},\varphi)\longmapsto (C_{\cal E},\varphi)$ identifies the set
of oriented pairs of type $(E,{\cal L})$ with the subspace $Z(j)$ of the
affine space
${\cal A}_\lambda\times A^0(E)$ which is cut-out by the integrability condition
$$j(C,\varphi):=(F^{02}_C,\bar\partial_C\varphi)=0
 $$

\begin{dt} A pair $(C,\varphi)\in  {\cal A}_\lambda\times A^0(E)$ will be called
\underbar{irreducible} if any $C$-parallel endomorphism $f\in A^0(su(E))$ with
$f(\varphi)=0$ vanishes.
\end{dt}
This notion of (ir)reducibility must not be confused with that one
introduced  in
section 1.2.3, which depends on the choice of an admissible pair. For instance,
irreducible pairs can be abelian.

The theorem above can now be reformulated as follows:
\begin{pr} An oriented pair $({\cal E},\varphi)$ of type $(E,{\cal L})$ is
polystable
if and only if the complex orbit ${\cal G}_0^\C\cdot (C_{\cal E},\varphi)\subset
Z(j)$ intersects the vanishing locus $Z(m)$ of the moment map $m$. $({\cal
E},\varphi)$ is stable if and only if it is  polystable and $(C_{\cal
E},\varphi)$ is
irreducible.
\end{pr}

It can be easily seen that the intersection $\left[{\cal G}_0^\C\cdot (C_{\cal
E},\varphi)\right]\cap Z(m)$ of a complex orbit with the vanishing locus of the
moment map  is either empty or coincides with a
\underbar{real}  orbit.  Moreover, using the proposition above one can
prove that the
set $Z(j)^{st}$ of stable oriented pairs is an \underbar{open} subset of the set
$Z(j)^{s}$ of simple oriented pairs. The quotient $\qmod{Z(j)^s}{{\cal
G}_0^{\C}}$ can
be  identified with the moduli space ${\cal M}^s(E,{\cal L})$ of simple
oriented pairs
of type $(E,{\cal L})$. The open subspace ${\cal M}^{st}(E,{\cal
L}):=\qmod{Z(j)^{st}}{{\cal G}_0^{\C}}\subset {\cal M}^s(E,{\cal L})$ will
be called
the moduli space of stable oriented pairs, and comes with a natural
structure of a
\underbar{Hausdorff} complex space.

The same methods as in [DK], [LT], [OT1] give finally the following
\begin{thry} The identification map $(C,\varphi)\longmapsto
(\bar\partial_C,\varphi)$ induces an isomorphism of real analytic spaces
$\qmod{Z(j,m)^{ir}}{{\cal G}_0}\textmap{\simeq}\qmod{Z(j)^{st}}{{\cal
G}_0^{\C}}={\cal M}^{st}(E,{\cal L})$, where
$Z(j,m)^{ir}$ denotes the locally closed subspace  consisting of irreducible
oriented pairs solving the equations $j(C,\varphi)=0$, $m(C,\varphi)=0$.
\end{thry}
\subsection{Decoupling the $PU(2)$-monopole equations}

Let $(X,g)$ be a K\"ahler surface and let $P^{\rm can}\map P_g$ be the
associated
\underbar{canonical} $Spin^c(4)$-\underbar{structure} whose spinor bundles are
$\Sigma^+=\Lambda^{00}\oplus\Lambda^{02}$,  $\Sigma^-=\Lambda^{01}$.  By
Propositions 1.1.11, 1.1.7 it follows that  the data of a
$Spin^{U(2)}(4)$-structure
in $(X,g)$ is equivalent to the data of a Hermitian 2-bundle $E$. The   bundles
associated with the
$Spin^{U(2)}(4)$-structure $\sigma:P^u\map P_g$ corresponding to $E$ are:
$$\det(P^u)=\det E\otimes K_X  ,\   \bar\delta(P^u)=\qmod{P_E}{S^1}  , $$ \
$$\Sigma^{\pm}(P^u)=\Sigma^{\pm}\otimes
E^{\vee}\otimes\det(P^u)=\Sigma^{\pm}\otimes E\otimes K_X\ ;\ \
\Sigma^{+}(P^u)=E
\otimes K_X\oplus E\ .
$$
 Suppose that $\det(P^u)\in NS(X)$ and fix an \underbar{integrable} connection
$a\in{\cal A}(\det(P^u))$. Denote by
$c\in {\cal A}(K_X)$ the Chern connection in $K_X$, by $\lambda:=a\otimes
\bar c$
the induced connection in $\det(E)=\det(P^u)\otimes \bar K_X$ and by ${\cal
L}$ the
corresponding holomorphic structure in this line bundle.  Identify the
affine space
${\cal A}(\bar\delta(P^u))$ with ${\cal A}_{\lambda\otimes c^{\otimes
2}}(E\otimes K_X)$  and the space of spinors
$A^0(\Sigma^+(P^u))$ with the direct sum $A^0(E \otimes K_X)\oplus
A^0(E)=A^0(E \otimes K_X)\oplus
A^{02}(E \otimes K_X)$.  The same computations  as in Proposition 4.1 [OT5]
gives
the   following {\it decoupling theorem}:
\begin{thry} A pair
$$(C,\varphi+\alpha)\in {\cal A}_{\lambda\otimes c^{\otimes 2}}(E\otimes
K_X)\times\left(A^0(E\otimes K_X)\oplus A^{02}(E \otimes K_X)\right)$$
solve the $PU(2)$-monopole equations $SW^\sigma_a$ if and only if the
connection
$C$ is integrable and one of the following conditions is fulfilled:
$$1)\ \alpha=0,\ \bar\partial_C\varphi=0\ \ \ and\ \ \ i\Lambda_g
F_C^0+\frac{1}{2}(\varphi\bar\varphi)_0=0\ ,
$$
$$\  2)\ \varphi=0,\  \partial_C\alpha=0\ \ and\ \ i\Lambda_g
F_C^0-\frac{1}{2}*(\alpha\wedge\bar\alpha)_0=0\ ,
$$
\end{thry}
Using Theorem 2.2.5 we get

\begin{re} The moduli space $({\cal M}^\sigma_a)_{\alpha=0}^{ir}$ of
irreducible
solutions of type 1) can be identified  with the moduli space ${\cal
M}^{st}(E\otimes
K_X,{\cal L}\otimes{\cal K}_X^{\otimes 2})$.

The moduli space $({\cal M}^\sigma_a)^{ir}_{\varphi=0}$ of
irreducible solutions of type 2) can be identified  with the moduli  space
${\cal
M}^{st}(E^{\vee},{\cal L}^{\vee})$ via the map $(C,\alpha)\longmapsto (\bar
C\otimes c,\bar \alpha)$.
\end{re}

Concluding, we get the following simple description of the moduli space ${\cal
M}^\sigma_a$ in terms of moduli spaces of stable oriented pairs.

\begin{co} Suppose that the $Spin^{U(2)}(4)$-structure $\sigma:P^u\map  P_g$ is
associated to the pair $(P^{\rm can}\map  P_g,E)$, where $P^{\rm can}\map  P_g$
is the canonical $Spin^c(4)$-structure of the K\"ahler surface $(X,g)$ and
$E$ is a
Hermitian rank 2 bundle. Let $a\in{\cal A}(\det(P^u))$ be an integrable
connection and
${\cal L}$ the holomorphic structure in $\det E=\det(P^u)\otimes
K_X^{\vee}$ defined by
$a$ and the Chern connection in $K_X$. Then the moduli space
${\cal M}^\sigma_a$ decomposes as a union of two Zariski closed subspaces
$${\cal M}^\sigma_a=({\cal M}^\sigma_a)_{\alpha=0}\union({\cal
M}^\sigma_a)_{\varphi=0}
$$
which intersect along the Donaldson moduli space ${\cal D}(\delta(P^u))
\subset{\cal M}^\sigma_a$ (see Remark 1.2.8).  There are canonical real analytic
isomorphisms
$$({\cal M}^\sigma_a)_{\alpha=0}^{ir}\simeq{\cal M}^{st} (E\otimes K_X,{\cal
L}\otimes{\cal K}_X^{\otimes 2})\ ,\ \ ({\cal
M}^\sigma_a)^{ir}_{\varphi=0}= {\cal
M}^{st}(E^{\vee},{\cal L}^{\vee})
$$
\end{co}

Using Remark 1.2.7, we recover the main result (Theorem 7.3) in [OT5]
stated for
quaternionic monopoles.
\vspace{5mm}\\
{\bf Example:} (R. Plantiko) On $\P^2$ endowed with the standard Fubini-Study
metric $g$ consider the $Spin^{U(2)}(4)$-structure $P^u\map  P_g$ with
$c_1(\det(P^u))=4$,
$p_1(\bar\delta(P^u))=-3$. It is easy to see that this
$Spin^{U(2)}(4)$-structure is
associated with  the pair $(P^{\rm can}\map P_g,E)$, where $E$ is a
$U(2)$-bundle with  $c_2(E)=13$, $c_1(E)=7$. Therefore $E\otimes K$ has
$c_1(E\otimes K)=1$, $c_2( E\otimes K)=1$.

Using Remark 2.1.4 it is easy to see
that any stable oriented pair $({\cal F},\varphi)$ of type $(E\otimes K,
{\cal O}(1))$
with
$\varphi\ne 0$ fits in an exact sequence of the form
$$
 0\map {\cal O} \textmap{\varphi} {\cal F}\map  {\cal O}(1)\otimes
J_{z_\varphi}\map 0\ ,$$
where $z_\varphi\in\P^2$, $c\in\C$ and ${\cal F}={\cal T}_{\P^2}(-1)$ is the
unique stable bundle with
$c_1=c_2=1$. Moreover,  two oriented pairs $({\cal F},\varphi)$,  $({\cal
F},\varphi')$ define the same point in the moduli space of stable oriented
pairs of
type $(E\otimes K,{\cal O}(1))$ if and only if $\varphi'=\pm \varphi$.
Therefore
$${\cal M}^{st}(E\otimes K,{\cal
O}(1))=\qmod{H^0({\cal F})}{\pm \id}\simeq \qmod{\C^3}{\pm\id}$$

Studying the local models of the moduli space one can check that the above
identification is a complex analytic isomorphism.

On the other hand every polystable oriented pair of type $ (E\otimes
K,{\cal O}(1))$
is stable and there is no polystable oriented pair of type $(E^{\vee},{\cal
O}(-7))$.
This shows that
$${\cal M}^\sigma_a\simeq\qmod{\C^3}{\pm\id} \
$$
for every integrable connection $a\in{\cal A}(\det(P^u))$. The quotient
$\qmod{\C^3}{\pm\id}$ has a  natural compactification   ${\cal
C}:=\qmod{\P^3}{\langle\iota\rangle}$, where $\iota$ is the involution
$$[x_0,x_1,x_2,x_3]\longmapsto [x_0,-x_1,-x_2,-x_3]\ .$$
${\cal C}$ can be identified with cone over the image of $\P^2$ under the
Veronese map
$v_2:\P^2\map \P^5$. This compactification coincides with the  {\it Uhlenbeck
compactification} of the moduli space [T1], [T2].
Let now $\sigma':P'^u\map  P_g$ be the $Spin^{U(2)}(4)$-structure in $\P^2$ with
$\det(P'^u)=\det(P^u)$, $p_1(\delta(P'^u))=+1$. It is easy to see by the
same method
  that ${\cal M}^{\sigma'}_a$ consists of only one point, which is the
{\it abelian}
solution associated with the {\it stable} oriented pair $({\cal O}\oplus{\cal
O}(1),\id_{\cal O})$.  Via the isomorphism explained in Proposition 1.2.15,
${\cal
M}^{\sigma'}_a$ corresponds to the
 moduli space of solutions of the (abelian) twisted Seiberg-Witten   equations
associated with the canonical $Spin^c(4)$-structure and the positive
chamber  (see
[OT6]).  Therefore
\begin{pr} The Uhlenbeck compactification of the moduli space ${\cal
M}^\sigma_a$   can be
identified  with the cone ${\cal C}$ over the image of $\P^2$ under the
Veronese
map $v_2$. The   vertex of the cone corresponds to the unique Donaldson
point. The
base  of the cone corresponds to the space
${\cal M}^{\sigma'}_a\times\P^2$ of ideal solutions concentrated in one
point.   The
moduli space ${\cal M}^{\sigma'}_a$ consists of only one abelian point.
\end{pr}

\newpage

\centerline{\large{\bf References}}
\vspace{6 mm}
\parindent 0 cm

[AHS] Atiyah M., Hitchin N. J., Singer I. M.: {\it Selfduality in
four-dimensional Riemannian geometry}, Proc. R. Lond. A. 362, 425-461 (1978)

[B] Bradlow, S. B.: {\it Special metrics and stability for holomorphic
bundles with global sections}, J. Diff. Geom. 33, 169-214 (1991)

[D] Donaldson, S.: {\it Anti-self-dual Yang-Mills connections over
complex algebraic surfaces and stable vector bundles}, Proc. London Math.
Soc. 3, 1-26 (1985)

[DK] Donaldson, S.; Kronheimer, P.B.: {\it The Geometry of
four-manifolds}, Oxford Science Publications   1990

[FU] Freed D. S. ;  Uhlenbeck, K.:
{\it Instantons and Four-Manifolds.}
Springer-Verlag 1984

[GS] Guillemin, V.; Sternberg, S.: {\it Birational equivalence in the
symplectic category},  Inv. math. 97, 485-522 (1989)

[HH] Hirzebruch, F.; Hopf,  H.: {\it Felder von Fl\"achenelementen in
4-dimensiona\-len 4-Mannigfaltigkeiten}, Math. Ann. 136 (1958)

[H] Hitchin, N.: {\it  Harmonic spinors}, Adv. in Math. 14, 1-55 (1974)

[HKLR] Hitchin, N.; Karlhede, A.; Lindstr\"om, U.; Ro\v cek, M.: {\it
Hyperk\"ahler
metrics and supersymmetry}, Commun.\ Math.\ Phys. (108), 535-589 (1987)

[K] Kobayashi, S.: {\it Differential geometry of complex vector bundles},
Princeton University Press  1987

[KM] Kronheimer, P.; Mrowka, T.: {\it The genus of embedded surfaces in
the projective plane}, Math. Res. Letters 1, 797-808 (1994)

[LL] Li, T.; Liu, A.: {\it General wall crossing formula}, Math. Res. Lett.
2,
  797-810  (1995).

[LM] Labastida, J. M. F.;  Marino, M.: {\it Non-abelian monopoles on
 four manifolds}, Preprint,
Departamento de Fisica de Particulas,   Santiago de Compostela, April
 (1995)

[La]  Larsen, R.: {\it Functional analysis, an introduction}, Marcel
Dekker, Inc., New York, 1973

[LMi] Lawson, H. B. Jr.; Michelson, M. L.: {\it Spin Geometry}, Princeton
University Press, New
Jersey, 1989

[LT] L\"ubke, M.; Teleman, A.: {\it The Kobayashi-Hitchin
correspondence},
 World Scientific Publishing Co.  1995

[M]  Miyajima, K.: {\it Kuranishi families of
vector bundles and algebraic description of
the moduli space of Einstein-Hermitian
connections},   Publ. R.I.M.S. Kyoto Univ.  25,
  301-320 (1989)

[MFK] Mumford, D,; Fogarty, J.; Kirwan, F.: {\it Geometric invariant
theory}, Springer Verlag,
1994

[OST]  Okonek, Ch.; Schmitt, A.; Teleman, A.: {\it Master spaces for stable
pairs}, Preprint,
alg-geom/9607015

[OT1]   Okonek, Ch.; Teleman, A.: {\it The Coupled Seiberg-Witten
Equations, Vortices, and Moduli Spaces of Stable Pairs},    Int. J. Math.
Vol. 6, No. 6, 893-910 (1995)

[OT2] Okonek, Ch.; Teleman, A.: {\it Les invariants de Seiberg-Witten
et la conjecture de Van De  Ven}, Comptes Rendus Acad. Sci. Paris, t.
321,  S\'erie I, 457-461 (1995)

[OT3] Okonek, Ch.; Teleman, A.: {\it Seiberg-Witten invariants and
rationality of complex surfaces}, Math. Z., to appear

[OT4] Okonek, Ch.; Teleman, A.: {\it Quaternionic monopoles},  Comptes
Rendus Acad. Sci. Paris, t. 321, S\'erie I, 601-606 (1995)

[OT5]  Ch, Okonek.;  Teleman, A.: {\it Quaternionic monopoles},
Commun.\ Math.\ Phys., Vol. 180, Nr. 2, 363-388, (1996)

[OT6] Ch, Okonek.;  Teleman, A.: {\it Seiberg-Witten invariants for
manifolds with
$b_+=1$,  and the universal wall crossing formula},
Int. J. Math., to appear

[PT1]  Pidstrigach, V.; Tyurin, A.: {\it Invariants of the smooth
structure of an algebraic surface arising from the Dirac operator},
Russian Acad. Izv. Math., Vol. 40, No. 2, 267-351 (1993)

[PT2]  Pidstrigach, V.; Tyurin, A.: {\it Localisation of the Donaldson
invariants along the
Seiberg-Witten classes}, Russian Acad. Izv. , to appear

[T1] Teleman, A. :{\it Non-abelian Seiberg-Witten theory},
Habilitationsschrift,
Universit\"at Z\"urich, 1996

[T2] Teleman, A. :{\it Moduli spaces of $PU(2)$-monopoles}, Preprint,
Universit\"at
Z\"urich, 1996

[W] Witten, E.: {\it Monopoles and four-manifolds}, Math.  Res.
Letters 1,  769-796  (1994)
\vspace{0.3cm}\\
Author's address : %
 Institut f\"ur Mathematik, Universit\"at Z\"urich,  Winterthu\-rerstr. 190,
CH-8057 Z\"urich,  {\bf  e-mail}:  teleman@math.unizh.ch\\

\hspace*{2.4cm} and Department of Mathematics, University of Bucharest.\\

\end{document}